\documentclass[12pt]{amsart}
\usepackage{amscd}

\newcommand{\C}{{\mathbb C}}
\newcommand{\Q}{{\mathbb Q}}
\newcommand{\M}{{\mathbb M}}
\newcommand{\Z}{{\mathbb Z}}
\newcommand{\R}{{\mathbb R}}
\newcommand{\F}{{\mathbb F}}
\newcommand{\PP}{{\mathbb P}}
\newcommand{\Qbar}{{\overline{\Q}}}
\newcommand{\Fbar}{{\overline{\F}}}
\newcommand{\Lbar}{{\overline{L}}}
\newcommand{\tors}{_{\text{tors}}}
\newcommand{\sep}{^{\text{sep}}}
\newcommand{\Char}{\operatorname{char}}
\newcommand{\rank}{\operatorname{rank}}
\newcommand{\End}{\operatorname{End}}
\newcommand{\Gal}{\operatorname{Gal}}
\newcommand{\GalQ}{{\Gal}(\Qbar/\Q)}

\newcommand{\Fq}{{{\mathbb F}_q}}
\newcommand{\defequal}{\stackrel{\text{def}}{=}}

\newcommand{\isom}{\cong}
\newcommand{\del}{\partial}
\newcommand{\calC}{{\mathcal C}}
\newcommand{\pp}{{\mathfrak p}}
\newcommand{\qq}{{\mathfrak q}}
\newcommand{\divisor}{\operatorname{div}}
\newcommand{\Div}{\operatorname{Div}}

\newcommand{\calF}{{\mathcal F}}
\newcommand{\calMt}{{{\mathcal M}_3}}
\newcommand{\calMp}{{{\mathcal M}_p}}
\newcommand{\calO}{{\mathcal O}}
\newcommand{\qp}{{\Q_p}}
\newcommand{\zt}{{\Z_3}}
\newcommand{\zp}{{\Z_p}}
\newcommand{\ft}{{\F_3}}

\newcommand{\cq}{{\calC(\Q)}}
\newcommand{\jq}{{J(\Q)}}
\newcommand{\jft}{{J(\F_3)}}
\newcommand{\jqt}{{J(\Q_3)}}
\newcommand{\jfp}{{J(\F_p)}}
\newcommand{\jqp}{{J(\Q_p)}}
\newcommand{\td}{{\widetilde D}}
\newcommand{\notdiv}{{\not\hskip-.2pt \vert\ }}

\newtheorem{lemma}{Lemma}

\newtheorem{prop}{Proposition}
\newtheorem{thm}{Theorem}

\theoremstyle{definition}
	
\newtheorem{conj}{Conjecture}
\newtheorem{question}{Question}	

\theoremstyle{remark}

\begin{document}

\title[Cycles of Quadratic Polynomials]{Cycles of Quadratic Polynomials and \\ Rational Points on a Genus~$2$ Curve}

\author{E.\ V.\ Flynn}
\address[Flynn]{Department of Pure Mathematics, University of Liverpool, P. O. Box 147, Liverpool L69 3BX, England}
\email{evflynn@liverpool.ac.uk}

\author{Bjorn Poonen}
\address[Poonen]{Mathematical Sciences Research Institute \\ Berkeley, CA 94720-5070, USA}
\email{poonen@msri.org}

\author[Schaefer]{Edward F.\ Schaefer}
\address{Santa Clara University \\ Santa Clara, CA 95053, USA}
\email{eschaefer@scuacc.scu.edu}

\thanks{The second author is supported by an NSF Mathematical Sciences Postdoctoral Research Fellowship.  Research at MSRI is supported in part by NSF grant DMS-9022140.  The third author is supported by an NSA Young Investigators Grant and a Paul Locatelli Junior Faculty Fellowship.}
\date{July 25, 1995}

\begin{abstract}
It has been conjectured that for $N$ sufficiently large, there are no quadratic polynomials in $\Q[z]$ with rational periodic points of period $N$.
Morton proved there were none with $N=4$, by showing that the genus~$2$ algebraic curve that classifies periodic points of period~4 is birational to $X_1(16)$, whose rational points had been previously computed.
We prove there are none with $N=5$.
Here the relevant curve has genus~$14$, but it has a genus~$2$ quotient, whose rational points we compute by performing a~$2$-descent on its Jacobian and applying a refinement of the method of Chabauty and Coleman.
We hope that our computation will serve as a model for others who need to compute rational points on hyperelliptic curves.
We also describe the three possible $\GalQ$-stable $5$-cycles, and show that there exist $\GalQ$-stable $N$-cycles for infinitely many $N$.
Furthermore, we answer a question of Morton by showing that the genus~$14$ curve and its quotient are not modular.
Finally, we mention some partial results for $N=6$.
\end{abstract}

\maketitle

\section{Introduction}
\label{intro}

Let $g(z) \in \Q(z)$ be a rational function of degree $d \ge 2$.
We consider $g$ as a map on $\PP^1(\C)$.
If $x \in \PP^1(\C)$ and the sequence
	$$x,g(x),g(g(x)),\ldots,g^{\circ n}(x),\ldots$$
is eventually periodic, then $x$ is called a {\em preperiodic point} for $g$.
If furthermore $g^{\circ n}(x) = x$, then $x$ is called a {\em periodic point} of $g$ of period $n$, and its orbit
	$$\{x,g(x),g(g(x)),\ldots,g^{\circ(n-1)}(x)\}$$
is called an {\em $n$-cycle} if $x$ does not actually have smaller period.
Northcott~\cite{northcott} proved in 1950 that for fixed $g$, there are only finitely many preperiodic points in $\PP^1(\Q)$.
Moreover, these can be computed effectively given $g$.
This theorem also holds over any fixed number field, and also for morphisms of $\PP^n$ of degree at least~2.
Since then, the theorem (in varying degrees of generality) has been rediscovered by many authors~\cite{narkiewicz}, \cite{lewis}, \cite{callsilverman}.

It is much more difficult to obtain {\em uniform} results for rational functions of a given degree.
Morton and Silverman~\cite{mortonsilverman} have proposed the following conjecture.
\begin{conj}
\label{bigconjecture}
Let $K/\Q$ be a number field of degree $D$, and let $\phi : \PP^N \rightarrow \PP^N$ be a morphism of degree $d \ge 2$ defined over $K$.  The number of $K$-rational preperiodic points of $\phi$ can be bounded in terms of $D$, $N$, and $d$ only.
\end{conj}

To demonstrate the strength of this conjecture, let us remark that the case $N=1$ and $d=4$ would imply the recently proved strong uniform boundedness conjecture for torsion of elliptic curves~\cite{merel}, namely that for any $D$ there exists $C>0$ such that for any elliptic curve $E$ over a number field $K$ of degree $D$ over $\Q$, $\# E(K)\tors < C$.
This is because torsion points of elliptic curves are exactly the preperiodic points of the multiplication-by-$2$ map, and their $x$-coordinates are preperiodic points for the degree~$4$ rational map that gives $x(2P)$ in terms of $x(P)$.
A similar conjecture for polynomials over $\Fq(T)$ and its finite extensions would imply the uniform boundedness conjecture for Drinfeld modules~\cite{poonen}, which is still open.

Even the simplest cases of the conjecture seem to be difficult.
Walde and Russo~\cite{walderusso} asked whether a quadratic polynomial in $\Q[z]$ can have rational points of period greater than~$3$, and this is not known.
Pairs consisting of a quadratic polynomial and a point of period $N$ are classified by an algebraic curve, which we denote $C_1(N)$.
For $N=1,2,3$, this curve is birational over $\Q$ to $\PP^1$, so it is easy to find a quadratic $g \in \Q[z]$ with a rational point of period $1$, $2$, or $3$.
Morton~\cite{morton4} proved that $C_1(4)$ is birational over $\Q$ to the modular curve $X_1(16)$, and used this to show that there are no quadratic polynomials in $\Q[z]$ with rational points of period~$4$.
Our main theorem is for the case $N=5$:

\begin{thm}
\label{fivetheorem}
There is no quadratic polynomial $g(z) \in \Q[z]$ with a rational point of exact period~5.
\end{thm}

The curve $C_1(5)$ has genus~$14$, so we study it via a quotient curve $\calC=C_0(5)$ of genus~$2$.
In Section~\ref{nonmodularity}, we will use the description of endomorphism rings of quotients of the Jacobian $J_1(N)$ of $X_1(N)$ to show that there is no surjective morphism of curves over $\C$ from $X_1(N)$ to $C_0(5)$ or $C_1(5)$, for any $N \ge 1$.
Because of this, finding the set of rational points will be more challenging than it was for $C_1(4)$.
To find all the rational points on $\calC$, we first put $\calC$ into hyperelliptic form, and then use a $2$-descent to compute the rank of its Jacobian, which turns out to be~$1$.
The $2$-descent is more difficult than the examples of descents for hyperelliptic curves worked out in the literature (\cite{flynndescent},\cite{GG},\cite{Sch}) in that $\calC$ has no Weierstrass points defined over $\Q$ or even a quadratic extension; in fact, the smallest field over which all the Weierstrass points are defined is the splitting field of a sextic with Galois group $S_6$, the worst possible case.
But because the rank is less than the genus, it is possible afterwards to apply the method of Chabauty and Coleman to bound the number of rational points on the curve.
Although Coleman's original method gives at best an upper bound of~$9$ for the number of rational points, our refinements of the method are strong enough to show that there are at most six rational points.
On the other hand, it is easy to list six rational points, so we know that we have found them all.

We will also list (in Table~\ref{cvalues}) all quadratic polynomials in $\Q[z]$ (up to linear conjugacy) with a $\GalQ$-stable $5$-cycle.
Each point in such a cycle generates a degree~$5$ cyclic extension of $\Q$, which we describe.
Also we prove that $\GalQ$-stable $N$-cycles exist for infinitely many $N$.

Finally, in Section~\ref{period6}, we describe the known $\GalQ$-stable $6$-cycles.
If, as we believe, these are all, then there is no quadratic polynomial $g(z) \in \Q[z]$ with a rational point of exact period~$6$.
This leads us to conjecture the following refinement of Conjecture~\ref{bigconjecture} for the case of quadratic polynomials over $\Q$.

\begin{conj}
\label{periodissmall}
If $N \ge 4$, then there is no quadratic polynomial $g(z) \in \Q[z]$ with a rational point of exact period~$N$.
\end{conj}

Throughout the paper, we use Mathematica (version 2.2) and the GP/PARI Calculator (version 1.39).
Version 1.39 of PARI assumes the Generalized Riemann Hypothesis for certain number field calculations, but Michel Olivier has kindly verified these particular calculations for us using a newer not yet released version that makes no such assumptions.

\section{Periodic points of quadratic polynomials}
\label{periodicpoints}

If $g(z) \in \Q[z]$ is any quadratic polynomial, then there exists a linear function $\ell(z) \in \Q[z]$ such that $\ell(g(\ell^{-1}(z)))$ is of the form $z^2+c$.
Therefore, for the sake of arithmetic dynamics, it will suffice to consider polynomials of the form $g(z)=z^2+c$.
If $z$ is periodic of exact period $N$ for $g$ (meaning that it is periodic of period $N$, but not periodic of period $n$ for any $n<N$), then $z$ satisfies the equation
\begin{equation}
\label{periodN}
	g^{\circ N}(z)-z=0.
\end{equation}
But (\ref{periodN}) is satisfied also by points of exact period $d$ for $d$ dividing $N$, so there is a factorization
	$$g^{\circ N}(z)-z = \prod_{d|N} \Phi_d(z,c)$$
where
\begin{equation}
\label{phidef}
	\Phi_d(z,c) = \prod_{m|d} (g^{\circ m}(z)-z)^{\mu(m)} \in \Z[z,c]
\end{equation}
is the polynomial whose roots $z$ for generic $c$ are the periodic points of exact period $d$.
(Here $\mu$ is the M\"{o}bius $\mu$-function.)
The $z$-degree of $\Phi_N(z,c)$ is
	$$\nu_2(N) \defequal \sum_{d|N} 2^d \mu(N/d).$$
By Theorem~1 in~\cite{bousch}, $\Phi_N(z,c)$ (where now $c$ also is considered to be an indeterminate) is irreducible in $\C[z,c]$, and hence
	$$\Phi_N(z,c)=0$$
defines a geometrically irreducible algebraic curve over $\Q$ in the $(z,c)$-plane.
Although the affine part of this curve is nonsingular (Proposition~1 of~\cite{bousch}), there is a singularity at infinity on its projective closure if $N>2$, so we let $C_1(N)$ denote the normalization, which is a nonsingular projective curve over $\Q$.
Every pair consisting of a polynomial $g(z)=z^2+c$ together with a rational point of exact period $N$ gives rise to a rational point on the affine part of $C_1(N)$.
The converse is true for almost all affine rational points, but there can be exceptions, as noted in Section~1 of~\cite{mortonsilverman2}, and these can be explained by assigning multiplicities to periodic points.
For example, $(z,c)=(-1/2,-3/4)$ is a point on $C_1(2)$, but $-1/2$ is actually a fixed point of $g(z)=z^2-3/4$ instead of a point of exact period~2.
(In fact, it seems likely that there are no other such examples for quadratic polynomials over $\Q$; this would follow from Conjecture~\ref{periodissmall}, for example.)

The curve $C_1(N)$ has an obvious automorphism $\sigma$ given in the $(z,c)$-plane by $(z,c) \mapsto (z^2+c,c)$.
(All we are saying here is that if $\alpha$ is a point of exact period $N$ for $g(z)=z^2+c$, then so is $g(\alpha)$.)
This automorphism generates a group $\langle \sigma \rangle$ of order $N$, and we let $C_0(N)$ be the quotient curve $C_1(N)/\langle \sigma \rangle$.
Then $C_0(N)$ is again a nonsingular projective curve over $\Q$, and its rational points correspond (with finitely many exceptions) to pairs consisting of a polynomial $g(z)=z^2+c$, $c \in \Q$, with a $\GalQ$-stable $N$-cycle.
For example, $C_0(4)$ has a rational point corresponding to $g(z)=z^2-31/48$ with the 4-cycle
$$\begin{CD}
	1/4 + \sqrt{-15}/6 	@>>>	 -1 + \sqrt{-15}/12	\\
	@AAA				@VVV			\\
	-1 - \sqrt{-15}/12	@<<<	1/4 - \sqrt{-15}/6
\end{CD}$$
(The notation is intended to remind the reader of the modular curves $X_0(N)$ and $X_1(N)$, which parameterize elliptic curves together with a cyclic subgroup of order $N$, or a point of order $N$, respectively.)
Because field automorphisms must preserve polynomial relations over $\Q$, the action of an automorphism in $\GalQ$ on a $\GalQ$-stable $N$-cycle is a rotation.
Thus we obtain a homomorphism $\GalQ \rightarrow \Z/N\Z$, and a point in such an $N$-cycle generates an abelian extension that is independent of which point was chosen, since any such point can be expressed as a polynomial over $\Q$ in any other.

Bousch~\cite{bousch} derived a formula for the genus of $C_1(N)$.
Later, Morton~\cite[Theorem C]{mortoncurves} generalized the formula to some other families of polynomials, and also derived a formula for the genus of $C_0(N)$, which is birational to his curve $\delta_N(x,c)=0$.
Here are the formulas, which are given in terms of $\nu(N) \defequal \nu_2(N)/2$:
	$$g(C_1(N)) = 1 + \left( \frac{N-3}{2} \right) \nu(N) - \frac{1}{2} \sum_{d|N, d \not= N} d \nu(d) \phi\left(\frac{N}{d}\right) ;$$
	$$g(C_0(N)) = 
		1+\left( \frac{1}{2} - \frac{3}{2N} \right) \nu(N) -
\frac{1}{2} \sum_{d|N, d \not= N} \nu(d) \phi\left(\frac{N}{d}\right)$$
if $N$ is odd; and
{\small
	$$g(C_0(N)) = 
		1+\left( \frac{1}{2} - \frac{3}{2N} \right) \nu(N) -
\frac{1}{2} \sum_{d|N, d \not= N} \nu(d) \phi\left(\frac{N}{d}\right) - \frac{1}{4N} \sum_{r|N,2|r,N/r \text{ odd}} \mu\left(\frac{N}{r} \right) 2^{r/2},
$$}%
if $N$ is even.
Table~\ref{genus} gives these values for $N \le 10$.

\begin{table}
\begin{center}
\begin{tabular}{|c||c|c|}
$N$	& $g(C_0(N))$	& $g(C_1(N))$	\\ \hline \hline
1	& 0		& 0		\\ \hline
2	& 0		& 0		\\ \hline
3	& 0		& 0		\\ \hline
4	& 0		& 2		\\ \hline
5	& 2		& 14		\\ \hline
6	& 4		& 34		\\ \hline
7	& 16		& 124		\\ \hline
8	& 32		& 285		\\ \hline
9	& 79		& 745		\\ \hline
10	& 162		& 1690		\\ \hline
\end{tabular}
\end{center}
\caption{Genus of $C_0(N)$ and $C_1(N)$ for $N \le 10$.}
\label{genus}
\end{table}

For $N=1$,~2, or~3, $C_1(N)$ is in fact birational over $\Q$ to $\PP^1$, so examples of quadratic polynomials in $\Q[x]$ with points of period~1,~2, or~3 exist in abundance.
These are classified explicitly in~\cite{walderusso}.
In~\cite{morton4}, it is proved that $C_1(4)$ is birational over $\Q$ to the curve
	$$v^2=u(u^2+1)(1+2u-u^2),$$
which also happens to be an equation for $X_1(16)$.
Although at first this may appear to be a surprising coincidence, we can give a partial explanation: the Jacobian of a genus~$2$ curve with an automorphism of order~$4$ defined over $\Q$ is automatically an abelian variety of $GL_2$-type, and hence conjecturally is a quotient of the Jacobian $J_1(N)$ of the modular curve $X_1(N)$ for some $N \ge 1$.
(See~\cite{ribetkorea}.)
It has been known since 1908 that (in modern terminology) no elliptic curve over $\Q$ has a rational point of order 16, so the only rational points of $X_1(16)$ are the rational cusps~\cite{levi}.
This fact is what enabled Morton~\cite{morton4} to prove that all rational points on $C_1(4)$ were at infinity.

Morton~\cite{morton4} asked whether $C_1(N)$ was modular also for $N>4$.
We will prove in Section~\ref{nonmodularity} that $C_0(5)$ and $C_1(5)$ are {\em not} modular.
The curve $C_1(5)$ is of genus~$14$ and is of degree~$30$ in the $(z,c)$-plane, so it is much too complicated to be studied directly.
Instead we will work with $\calC \defequal C_0(5)$, which has genus~2.
Of course, every rational point of $C_1(5)$ maps to a rational point of $\calC$.

Before proceeding with the calculation of the rational points of $\calC$, let us show that the ``affine part'' of $C_0(N)$ has rational points for infinitely many $N$.
This contrasts with the modular curve situation, since for $N>163$, the only rational points of $X_0(N)$ are the rational cusps.
(The result for $X_0(N)$ involved many cases, which were worked out by several different authors.  See~\cite{kenku} for a brief summary.)

\begin{thm}
\label{easy}
There are infinitely many $N$ for which there exists a quadratic polynomial $g(z) \in \Q[z]$ with a $\GalQ$-stable $N$-cycle.
\end{thm}

\begin{proof}
For each $k \ge 1$, the image of $2$ is a generator of $(\Z/3^k\Z)^\ast$.
Then under the map $g(z)=z^2$, the orbit of a primitive $3^k$-th root of unity $\zeta$ is a $(2\cdot 3^{k-1})$-cycle consisting of all primitive $3^k$-th roots of unity, which is clearly $\GalQ$-stable.
(A similar argument could be used with $g(z)=z^2-2$ and $\zeta+\zeta^{-1}$.)
\end{proof}

Although the proof was disappointingly simple, it does raise an interesting question.
\begin{question}
Is it true that for sufficiently large $N$, if $g(z)=z^2+c$ has a $\GalQ$-stable $N$-cycle, then $c=0$ or $c=-2$?
\end{question}
For many $N$ (for example, $N=7$), not even $z^2$ and $z^2-2$ have $\GalQ$-stable $N$-cycles.
More precisely, it is easy to show that $z^2$ has a $\GalQ$-stable $N$-cycle if and only if $N=\phi(n)$ where $n$ is a positive integer for which the image of $2$ is a generator of $(\Z/n\Z)^\ast$ (which forces $n$ to an odd prime power).
Similarly, $z^2-2$ has a $\GalQ$-stable $N$-cycle if and only if $N=\phi(n)/2$ where the image of $2$ generates $(\Z/n\Z)^\ast/\langle -1 \rangle$ (which forces $n$ to be the product of at most two odd prime powers).

\section{A hyperelliptic form of $\calC$}
\label{hyperelliptic}

Because $\calC$ has genus~2, it is hyperelliptic, and more specifically is birational to a curve $\calC$ of the form $y^2 = f(x)$, where $f(x) \in \Q[x]$ is of degree~5 or~6 and has distinct roots.
For the future calculations, it will be necessary to find $f(x)$ explicitly.
This will be the concern of this section.

Following Morton~\cite{morton4}, we define the {\em trace} of an $N$-cycle in $\C$ of $g(z)=z^2+c$ to be the sum of the elements in the cycle.
Then we let $\tau_N(z,c) \in \Z[z,c]$ be the polynomial whose roots for generic $c$ are the traces of all the $N$-cycles.
The curve $\tau_N(z,c)=0$ is birational over $\Q$ to $C_0(N)$.
(See~\cite{mortoncurves}.)
In~\cite{morton4}, Morton also gives an efficient method for computing $\tau_N(z,c)$ for small $N$.

We will start with his result for $N=5$:
\begin{equation*}
\begin{split}
	\tau_5(z,c)	&=	\left( 32 + 28c + 40{c^2} + 9{c^3} \right) + 
				\left( 36 - 24c + 17{c^2} \right) z + 
				\left( 44 + 19c + 19{c^2} \right) {z^2} \\
			&\quad	+ \left( 11 + 18c \right) {z^3} + 
				\left( 3 + 11c \right) {z^4} + {z^5} + {z^6}.
\end{split}
\end{equation*}
Solving the system
	$$\tau_5 = \del \tau_5/\del z = \del \tau_5/\del c = 0,$$
we find that the only singularity of the curve $\tau_5(z,c)=0$ in the affine $(z,c)$-plane is $(-1,-4/3)$, which is a node.
Therefore we substitute $z=r-1$ and $c=s-4/3$ and clear denominators to obtain a new model with the node at $(0,0)$:
{\footnotesize
$$	238{r^2} + 213{r^3} - 15{r^4} - 45{r^5} + 9{r^6} + 36rs - 
   177{r^2}s - 234{r^3}s + 99{r^4}s + 54{s^2} - 189r{s^2} + 
   171{r^2}{s^2} + 81{s^3} = 0.$$}%
Next we blow up the node by substituting $s=rt$, and dividing by $r^2$:
{\small
$$238 + 213r - 15{r^2} - 45{r^3} + 9{r^4} + 36t - 177rt - 234{r^2}t + 99{r^3}t + 54{t^2} - 189r{t^2} + 171{r^2}{t^2} + 81r{t^3} = 0.$$}%
The curve now has no affine singularities, but there must be a singularity at infinity, because a nonsingular plane curve cannot have genus~2.
A calculation shows that there is a singularity at infinity on the line $r+t=0$, which we move to an axis by setting $r=q-t$:
$$238 + 213q - 15{q^2} - 45{q^3} + 9{q^4} - 177t - 147qt - 99{q^2}t + 63{q^3}t + 216{t^2} + 144q{t^2} - 72{q^2}{t^2} = 0.$$
Now the left hand side is a quadratic in $t$, so the curve is birational to
	$$p^2 = -174303 - 269082q + 15471{q^2} + 115668{q^3} + 5103{q^4} - 
   30618{q^5} + 6561{q^6},$$
where the right hand side is the discriminant of that quadratic.
Although this is a hyperelliptic form, it is to our advantage to simplify as much as possible before continuing.
We substitute $p=192y$ and $q=-1-4x/3$, and cancel $192^2=36864$ from both sides to obtain
\begin{equation}
\label{curveequation}
	\calC: \quad y^2 = x^6 + 8 x^5 + 22 x^4 + 22 x^3 + 5 x^2 + 6 x + 1.
\end{equation}
Let $f(x)$ be the sextic on the right hand side.
Since $f(x)$ has no rational roots, the curve $\calC$ is not birational over $\Q$ to a curve of the form $y^2 = h(x)$ with $\deg h(x) = 5$.

\section{Six rational points on $\calC$}
\label{sixpoints}

There are a few easy to find rational points on $\calC$.
First of all, $f(0)=f(-3)=1$, so we find four affine points: $(0,1)$, $(0,-1)$, $(-3,1)$, and $(-3,-1)$.
Also, since $\deg f$ is even, $\calC$ has two points at infinity.
Since the leading coefficient of $f(x)$ is a square in $\Q$, these points are rational.
(See~\cite[p.\ 50]{CA2}.)
The rational function $y/x^3$ takes values~1 and~$-1$ at these two points, which we call $\infty^+$ and $\infty^-$, respectively.

We will eventually show that these six points are the only rational points on $\calC$.
For now, we will describe the $5$-cycles of quadratic polynomials to which they correspond.
By tracing back through the substitutions of Section~\ref{hyperelliptic}, we obtain two equivalent formulas for $c$ in terms of the rational functions $x$ and $y$ on $\calC$:
$$c=\frac{P_0(x)+P_1(x)y}  {8{x^2}{{\left( 3 + x \right) }^2}}
=
\frac{64 + 110x + 325{x^2} + 452{x^3} + 271{x^4} + 74{x^5} + 8{x^6}}{2(P_0(x)-P_1(x)y)}
,$$
where
\begin{align*}
	P_0(x)	&= -9 - 24x - 95{x^2} - 104{x^3} - 46{x^4} - 10{x^5} - {x^6}	\\
	P_1(x)	&= -9 + 3x + 6{x^2} + {x^3}.
\end{align*}

The second formula is determinate (i.e., the numerator and denominator do not both vanish) at the four affine rational points, and this gives the $c$-values shown in Table~\ref{cvalues}.
At $\infty^+$, we have the formal expansion
	$$y = x^3 + 4 x^2 + 3x - 1 + 2 x^{-1} + \cdots.$$
Substituting this into the first formula, we see that
	$$c = x^2/4 + \text{(lower order terms)}$$
so $c$ has a pole at $\infty^+$.
Similarly, at $\infty^-$, we have
	$$y = -(x^3 + 4 x^2 + 3x - 1 + 2 x^{-1} + \cdots),$$
and substitution into the second formula shows that $c=-2$ there.

\begin{table}
\begin{center}
\begin{tabular}{|c||c|c|c|}
Point		& $c$		& Conductor $n$	& $\Gal(\Q(\zeta_n)/K)$	\\ \hline \hline
$(0,1)$		& $\infty$	&		&		\\ \hline
$(0,-1)$	& $-16/9$	& $41$		& $\langle 3 \rangle \subset (\Z/41\Z)^\ast$	\\ \hline
$(-3,1)$	& $-64/9$	& $275=5^2 \cdot 11$	& $\langle -1,3 \rangle \subset (\Z/275\Z)^\ast$	\\ \hline
$(-3,-1)$	& $\infty$	&		&		\\ \hline
$\infty^+$	& $\infty$	&		&		\\ \hline
$\infty^-$	& $-2$		& $11$		& $\langle -1 \rangle \subset (\Z/11\Z)^\ast$	\\ \hline
\end{tabular}
\end{center}
\caption{The six rational points of $\calC$.}
\label{cvalues}
\end{table}

For the three values $c=-2, -16/9, -64/9$, we know there is a $\GalQ$-stable 5-cycle of $g(z)=z^2+c$.
The action of $\GalQ$ on the cycle can only be a rotation, so the points of the cycle generate an abelian extension $K$ of $\Q$, whose Galois group is a subgroup of $\Z/5\Z$.
In Table~\ref{cvalues}, we will describe $K$ in each case by giving its conductor (the smallest $n$ for which $K$ is contained in the $n$-th cyclotomic field $\Q(\zeta_n)$) and the subgroup $\Gal(\Q(\zeta_n)/K)$ of $\Gal(\Q(\zeta_n)/\Q) \isom (\Z/n\Z)^\ast$ it corresponds to under Galois theory.
The quintic polynomial whose roots are the points of the cycle is a factor of $\Phi_5[z,c]$.
A computation shows that for each of the three values of $c$ above, there is a unique quintic factor in $\Q[z]$, and none of smaller degree, so already we know that the 5-cycles in question are not defined pointwise over $\Q$, and that in each case $K$ is a degree~5 cyclic extension of $\Q$.

For $c=-2$, PARI tells us that the field $K$, which is generated by a root of this quintic, has discriminant $11^4$, so the conductor of $K$ must be a power of $11$.
Since $(\Z/11^k\Z)^\ast$ is cyclic, $\Q(\zeta_{11^k})$ has a unique quintic subfield, namely the totally real subfield of $\Q(\zeta_{11})$.
Thus the conductor of $K$ equals 11, and under Galois theory $K$ corresponds to the subgroup $\langle -1 \rangle$ of $(\Z/11\Z)^\ast$.
This is easy to explain: the 5-cycle of $z^2-2$ consists of all conjugates of $\zeta_{11}+\zeta_{11}^{-1}$.

For $c=-16/9$, $K$ has discriminant $41^4$, so a similar argument as for $c=-2$ shows that $K$ is the unique quintic subfield of $\Q(\zeta_{41})$.
Thus $K$ has conductor 41, and corresponds to the unique subgroup of $(\Z/41\Z)^\ast$ of index 5, which is generated by the image of $3$.

For $c=-64/9$, $K$ has discriminant $5^8 \cdot 11^4$, so the conductor of $K$ is of the form $n=5^k \cdot 11^l$.
By Hensel's Lemma, every element of $(\Z/n\Z)^\ast$ congruent to $1$ modulo $275=5^2 \cdot 11$ is a 5-th power in $(\Z/n\Z)^\ast$, and hence is in $H \defequal \Gal(\Q(\zeta_n)/K)$, which has index~5 in $(\Z/n\Z)^\ast$.
Thus $n$ divides $275$.
PARI tells us that the prime 3 splits completely in $K$, so the Frobenius element at 3 acts trivially on $K$, and the image of 3 lies in $H$.
Also, the image of $-1$ lies in $H$, since $H$ has odd index.
But the subgroup generated by $-1$ and $3$ in $(\Z/275\Z)^\ast$ has index 5, so the images of $-1$ and $3$ in $(\Z/n\Z)^\ast$ generate $H$.
Finally, this subgroup of $(\Z/275\Z)^\ast$ does not come from a subgroup of $(\Z/55\Z)^\ast$, so the conductor is actually 275.

\section{Generalities on $2$-descents on Jacobians of hyperelliptic curves}
\label{generalities}

This section outlines and elaborates upon the descent method described in~\cite{CA2} for Jacobians of genus~$2$ curves over $\Q$.
(See also~\cite{flynndescent},~\cite{GG} and~\cite{Sch}.)
Later, in Section~\ref{descent}, we will apply the results of this section to show that the Mordell-Weil rank of the Jacobian of our curve $\calC$ is exactly~1.
We hope that the separation of the general method from the application will be useful for others who need to do $2$-descents on hyperelliptic curves.

Let $C$ be a hyperelliptic curve over $\Q$ of genus $g \ge 2$.
Then $C$ has a (singular) plane model $y^2=f(x)$, with $f(x) \in \Z[x]$ a separable polynomial of even degree $d=2g+2$.
Let $J$ be the Jacobian of $C$, which is an abelian variety over $\Q$.
We will assume $C(\Q)$ is nonempty, so that $\Div^0(C)(K)$ maps onto $J(K)$ for any field extension $K$ of $\Q$.
(Actually, when $g=2$, the latter is true automatically, even when $C(\Q)$ is empty.)
Without this assumption, the map $(x-T)$ below could be defined only as a map on $\Div^0(C)(K)$.
We will call a degree~$0$ divisor of $\calC$ defined over $\Q$ a {\em good divisor} if its support does not include $\infty^{+}$, $\infty^{-}$ or points with $y$-coordinate 0.

\begin{prop}
\label{good}
Every divisor class of $J(K)$ contains a good divisor.
\end{prop}

\begin{proof}
Since $C$ has a $K$-rational point, every
$K$-rational divisor class contains a $K$-rational divisor
(see~\cite[p.\ 168]{MI}). Every $K$-rational divisor has a linearly
equivalent $K$-rational divisor whose support avoids any given finite
set of points (see~\cite[p.\ 166]{LA}).
\end{proof}

For a good divisor $D=\sum n_{P}P$, we define
	$$(x-T)(D)=\prod_P (x_P-T)^{n_{P}}\in L_{K}^{\ast}.$$
(We use the notation $P=(x_P,y_P)$.)
For any field $K$ of characteristic $0$, define $L_{K}=K[T]/(f(T))$.

\begin{prop}
\label{kernorm}
The map $(x-T)$ is a well-defined map from $J(K)$ to the kernel of
the norm from 
$L_{K}^{\ast}/L_{K}^{\ast 2}{K}^{\ast}$ to $K^{\ast}/K^{\ast 2}$.
\end{prop}

\begin{proof}
Let $\alpha_1,\ldots,\alpha_d$ be the zeros of $f(x)$ in $\overline{K}$.
We can define
	$$\Lbar_K=\overline{K}[T]/(f(T))\cong \overline{K}[T]/(T-\alpha_1)
	\times \ldots \times \overline{K}[T]/(T-\alpha_{d})\cong
	\overline{K} \times\ldots\times\overline{K}$$
	$${\rm by}\;\; T\mapsto (\alpha_{1},\ldots ,\alpha_{d}).$$
Let $\Gal(\overline{K}/K)$ act trivially on $T$; this
makes $\Lbar_K$ a $\Gal(\overline{K}/K)$-module
and $L_{K}$ is the set of
$\Gal(\overline{K}/K)$-invariants.
Then we can consider $(x-T)$ to be a
$\Gal(\overline{K}/K)$-invariant $d$-tuple of functions
$((x-\alpha_{1}),\ldots ,(x-\alpha_{d}))$ whose divisors are
$(2(\alpha_{1},0)-\infty^{+}-\infty^{-}, \ldots ,2(\alpha_{d},0)-
\infty^{+}-\infty^{-})$.
We denote this $d$-tuple of divisors by $2(T,0)-\infty^{+}
-\infty^{-}$.

To show that $(x-T)$ is a well-defined map from $J(K)$ to
$L_{K}^{\ast}/L_{K}^{\ast 2}K^{\ast}$, we first note from
Proposition~\ref{good} that every element
of $J(K)$ contains a good divisor.
Let $D_1$ and $D_2$ be two good divisors that are linearly equivalent.
Then there is a $K$-defined function $h$ with $D_1-D_2=
\divisor h$. We have
the following equalities of $d$-tuples:
\[
(x-T)(D_1-D_2)=(x-T)(\divisor h)=h(\divisor (x-T))=
h(2(T,0)-\infty^+ - \infty^{-})=\]\[
h((T,0))^{2}/h(\infty^{+}) h(\infty^{-})\in L_{K}^{\ast 2}K^{\ast}.\]

Now let us show that the image of $(x-T)$ is contained in the kernel of the
norm to $K^{\ast}/K^{\ast 2}$.
Let $D=\sum n_{P}P$ be a good divisor. 
If $c$ is the leading coefficient of $f(x)$, then
	$$N_{L_{K}/K}((x-T)(D))
	=\prod_{P} \prod_{j=1}^{d}(x_P-\alpha_{j})^{n_{P}}
	=\prod_{P} (y_P^2/c)^{n_{P}}
	=\left( \prod_{P} y_P^{n_{P}} \right)^{2}\in K^{\ast 2}.$$
\end{proof}

Let $L=L_\Q=\Q[T]/((f(T)) \isom \prod_{i=1}^{r} L_{i}$, where the $L_i$ are fields corresponding to the irreducible factors of $f(x)$.
Let $S$ be a finite set of primes of $\Q$ containing the primes $2$, $\infty$, and all primes dividing the discriminant of $f(x)$.
(In particular, $S$ contains all primes dividing the leading coefficient of $f(x)$.)
Suppose $l \in L^\ast$ maps to $l_i$ in $L_i^\ast$.
Then we say that $l$ is {\em unramified outside $S$} if for each $i$,
the field extension
$L_{i}(\sqrt{l_{i}})/L_{i}$ is unramified outside of
primes lying over primes of $S$.
This property of $l$ depends only on the image of $l$ in $L^\ast/L^{\ast 2}$,
and it is easy to see that the subset $G$ of elements of $L^\ast/L^{\ast 2}$
which are unramified outside $S$ is a subgroup.
Let $G'$ be the image of $G$ in $L^{\ast}/L^{\ast 2}\Q^{\ast}$, and
let $H$ be the kernel of the norm from $G'$ to $\Q^{\ast}/\Q^{\ast 2}$.
 
\begin{prop}
\label{landsinh}
The image of the map $(x-T)$ on $J(\Q)$ is contained in the subgroup $H$ of $L^{\ast}/L^{\ast 2}\Q^{\ast}$.
\end{prop}

\begin{proof}
By Proposition~\ref{kernorm}, the image of $(x-T)$ is contained in the kernel of the norm to $\Q^\ast/\Q^{\ast 2}$.
So it suffices to show that the image of $(x-T)$ on any good divisor $D = \sum_P n_P P$ is contained in $G'$.

For each $p \not \in S$, fix an embedding $\Qbar \rightarrow \Qbar_p$.
Let $v$ be the additive $p$-adic valuation on $\Qbar_p$ with $v(p)=1$.
Since $D$ is $\GalQ$-stable, $\prod_{v(x_P)<0} x_P^{n_P}$ is fixed by the inertia group of $\Gal(\Qbar_p/\Q_p)$ and hence its valuation is an integer $a_p$.
Moreover since the embedding $\Qbar \rightarrow \Qbar_p$ is unique up to the action of $\GalQ$ on the left, $a_p$ is independent of the embedding.

Let $m=\prod_{p \not\in S} p^{a_p} \in \Q^\ast$.
We claim that $m^{-1} (x-T)(D) \in L^\ast/L^{\ast 2}$ is unramified outside $S$ (i.e., is in $G$), or what is the same thing, that for any $p \not\in S$ and any ring homomorphism $\iota: \Qbar[T]/(f(T)) \rightarrow \Qbar_p$, $v(m^{-1}(x-T)(D))$ is an even integer.
(We extend $v$ to $\Qbar[T]/(f(T))$ by applying $\iota$ when necessary.)
Let $\alpha_1,\ldots,\alpha_d$ be the zeros of $f(x)$ in $\Qbar_p$, and without loss of generality assume $\iota(T)=\alpha_1$.
If $v(x_P-T)>0$, then $v(x_P-\alpha_i)=0$ for $2 \le i \le d$, since the $\alpha_i$ lie in distinct residue classes of the ring of integers of $\Qbar_p$.
In this case,
	$$v(x_P-T) = v \left(\prod_{i=1}^d (x_P-\alpha_i) \right) = v(y_P^2/c) = 2 v(y_P),$$
where $c$ is the leading coefficient of $f(x)$, which by assumption is an $S$-unit.
Hence
\begin{align}
	v((x-T)(D))	&= \sum_{v(x_P-T)>0} v((x_P-T)^{n_P}) + \sum_{v(x_P-T)<0} v((x_P-T)^{n_P})	\\
\label{endformula}	&= 2 v \left( \prod_{v(x_P-T))>0} y_P^{n_P} \right) + v(m)
\end{align}
since $v(x_P-T)=v(x_P)$ when either is negative.
The product in~(\ref{endformula}) is again stable under the inertia group of $\Gal(\Qbar_p/\Q_p)$, so its valuation is an integer.
Thus $v(m^{-1}(x-T)(D))$ is an even integer.
\end{proof}

Let $L_p = L_{\Q_p} = \Q_p[T]/(f(T))$.
We have a commutative diagram
\begin{equation}
\label{cd}
\begin{CD}
	0 @>>> J(\Q)/\ker(x-T)	@>{x-T}>>	L^\ast/L^{\ast 2} \Q^\ast	\\
	&&	@VVV					@VVV	\\
	0 @>>> \prod_{p \in S} J(\Q_p)/\ker(x-T)	@>{x-T}>>	\prod_{p \in S} L_p^\ast/L_p^{\ast 2} \Q_p^\ast.
\end{CD}
\end{equation}
{}From this diagram and Proposition~\ref{landsinh},
we deduce that $x-T$ maps $J(\Q)/\ker(x-T)$ injectively into the subgroup $H'$ of elements of $H$ that for each $p \in S$ map in $L_p^\ast/L_p^{\ast 2} \Q_p^\ast$ into the image of $J(\Q_p)$ under $x-T$.
The latter is something that can be calculated, and this will give bounds on the size of $J(\Q)/\ker(x-T)$.

In order to convert these bounds into bounds on the size of $J(\Q)/2J(\Q)$, which will let us bound the rank of $J(\Q)$, we need to know how $\ker(x-T)$ compares with $2J(\Q)$.
Since $(x-T)$ maps $J(\Q)$ into an elementary $2$-group, clearly $2J(\Q) \subseteq \ker(x-T)$.
We will describe the difference between these two groups in Proposition~\ref{halves} below, for the genus $2$ case.
The result is stated over arbitrary fields of characteristic not equal to $2$, since we will need it for the completions of $\Q$ as well as for $\Q$ itself.
We will make use of the following well known consequence of the Riemann-Roch theorem.

\begin{prop}
\label{RR}
Suppose $\deg f(x)=6$, so the genus of $C$ is $2$.
Then any divisor class in $J(K)$ may be represented by a divisor of the form
$P_1 + P_2 - \infty^+ - \infty^-$ where
either $P_1,P_2\in C(K)$ or
$P_1,P_2 \in C (K')$, with $[K':K ] = 2$
and $P_1,P_2$ conjugate over $K$. 
This representation is unique (up to interchanging $P_1$ and $P_2$),
except for the group identity $\calO$ of $J(K)$, which can be represented by any divisor of the form $(x,y)+(x,-y)-\infty^{+}-\infty^{-}$ or $\infty^+ + \infty^-  - \infty^+ - \infty^- $.
\end{prop}

\begin{prop}
\label{halves}
Suppose that $f(x) \in K[x]$ is a separable sextic polynomial over a field $K$ with $\Char(K)\neq 2$, and that the
genus $2$ curve $C:y^{2}=f(x)$ has a point $P$ defined over $K$.
Let $J$ be the Jacobian of $C$.
Then the index of $2J(K)$ in $\ker(x-T)$ is 
\begin{quote}
	1 $\;\;$ \parbox[t]{5 in}{if $f(x)$ has a zero in $K$, or if there is some $\Gal(K\sep/K)$-stable partition of the six zeros into two indistiguished 3-element subsets $\{ \{ \alpha_{1}, \alpha_{2}, \alpha_{3} \} , \{\alpha_{4}, \alpha_{5}, \alpha_{6} \} \}$} \\
	2 $\;\;$ \hbox{otherwise.}
\end{quote}
\end{prop}

\begin{proof}
The index of $2J(K)$ in $\ker(x-T)$ is~$1$ or~$2$
and $\ker(x-T)/2J(K)$ is generated
by $[2P-\infty^{+}-\infty^{-}]$
(see~\cite[lemma 5.2,theorem 5.3]{CA2}).
So the index is~$1$ exactly when
$[2P-\infty^{+}-\infty^{-} ]$ is in $2J(K)$.
Now $[2P-\infty^{+}-\infty^{-} ]$ is in $2J(K)$
if and only if one of the 16 divisor classes with double
$[2P-\infty^{+}-\infty^{-} ]$ is in $J(K)$.

We now find these $16$ divisor classes.
Let $\alpha_{1},\ldots ,\alpha_{6}$ be the roots of $f(x)$ in some algebraic closure.
We will use repeatedly and without further mention the fact that the divisors $2(\alpha_i,0)$ and $\infty^+ + \infty^-$ are linearly equivalent.
Since
	$$2[P + (\alpha_{1},0)- \infty^+ - \infty^{-} ] = [2P-\infty^+ -\infty^-],$$
the 16 halves of $[2P- \infty^{+}-\infty^{-} ]$ can be obtained by
adding $[P+(\alpha_1,0)-\infty^+ -\infty^{-} ]$ to each of the 16 elements
of $J[2]$.
By Proposition~\ref{RR},
the 15 divisor classes $[ (\alpha_{i},0)+(\alpha_{j},0)-\infty^{+}-
\infty^{-} ]$ with
$i < j$ are distinct, and each has order 2.
Thus the 16 halves of $[2P - \infty^{+}-\infty^{-} ]$ are the~$6$ divisor
classes of the form
	$$[P+2(\alpha_1,0)+(\alpha_i,0)-2\infty^+-2\infty^-] = [P+ (\alpha_{i},0)-\infty^+ -\infty^{-}]$$
and the 10 divisor classes of the form
	$$[P+ (\alpha_{1},0)+(\alpha_{j},0)+(\alpha_{k},0)-2\infty^{+}-2\infty^{-} ]$$
with $1 < j < k$.

The action of $\Gal(K\sep/K)$ on the first~$6$ halves is the same as the action on the roots $\alpha_1,\ldots,\alpha_6$.
To deduce the action on the other 10 halves, note that if $1<j<k$ and $l,m,n$ are the other three possible indices, then
{\small
	$$[P+ (\alpha_l,0)+(\alpha_m,0)+(\alpha_n,0)-2\infty^{+}-2\infty^{-} ] = [P+ (\alpha_1,0)+(\alpha_{j},0)+(\alpha_{k},0)-2\infty^{+}-2\infty^{-} ]$$}%
because the difference of the two divisors is $\divisor((x-\alpha_1)(x-\alpha_j)(x-\alpha_k)/y)$.
Hence the action of $\Gal(K\sep/K)$ on these 10 halves is the same as the action on the 10 partitions of the six roots into two indistinguished 3-element subsets.

Thus the conditions given in the proposition are necessary and sufficient for $[2P-\infty^+-\infty^-]$ to be in $2J(K)$.
By our earlier remarks, this completes the proof.
\end{proof}

We conclude this section with a few remarks on computing the function $(x-T)$.
Although $P+Q-\infty^+-\infty^-$ is not a good divisor, the image of $(x-T)$ on its divisor class can be found in terms of $P$ and $Q$.
This is described in~\cite[p.\ 50]{CA2}. As an example, if $P$ and $Q$ are
both affine and have nonzero $y$-coordinates, then the image of 
$[ P + Q -\infty^+ -\infty^- ]$ is $(x_P-T)(x_Q-T)$.
In addition, the image of $[P+\infty^{\pm} -\infty^{+}-\infty^{-}]$
is $(x_P-T)$.

\section{Facts about the number field $L=\Q[T]/(f(T))$}
\label{numberfield}

{}From now on, we specialize to our curve $\calC$, for which
	$$f(x)=x^6 + 8 x^5 + 22 x^4 + 22 x^3 + 5 x^2 + 6 x + 1.$$
Let $L=\Q[T]/(f(T))$.
(We will abuse notation by writing $T$ for its image in $L$.)
In this section we will record some data on $L$ obtained from PARI, to be used later, mainly for the 2-descent.
The polynomial $f(x)$ is irreducible over $\Q$, so $L$ is a number field.
The Galois group of the normal closure $M$ of $L$ is the full symmetric group $S_6$.
The class number of $L$ is~1.
(This can be verified without using PARI, without too much difficulty, since the Minkowski bound is only about $12.2$.)
Two of the six zeros of $f(x)$ are real, so the unit group $U$ has rank~3.
The torsion of the unit group is only $\{\pm 1\}$, and the quotient $U/\{\pm 1\}$ is generated by the elements $u_1,u_2,u_3$ listed in Table~\ref{elements}.
The discriminant of $f$ is $2^{12} \cdot 3701$, and the prime factorizations of the ramified prime ideals $(2)$ and $(3701)$ in $L$ are $(\alpha)^2$ and $(\beta_1)(\beta_2)^2(\beta_3)$, respectively, where $\alpha,\beta_1,\beta_2,\beta_3$ are defined as in Table~\ref{elements}.
The factorization of $2$ and $3701$ into irreducible {\em elements} of $L$ will be given in Table~\ref{badprimes}.

\begin{table}
\begin{center}
\begin{tabular}{|c|c|c|}
Element		& Definition				& Norm \\ \hline \hline
$u_1$		& $(T^3+4T^2+3T-1)/2$			& $1$	\\ \hline
$u_2$		& $(T^4+5T^3+7T^2+2T+1)/2$		& $1$	\\ \hline
$u_3$		& $(T^4+6T^3+11T^2+5T)/2$		& $-1$	\\ \hline
$-1$		& $-1$					& $1$	\\ \hline
$\alpha$	& $(T^5+8T^4+22T^3+23T^2+7T+5)/2$	& $2^3$	\\ \hline
$\beta_1$	& $(-T^5-5T^4-5T^3+2T^2-3T+6)/2$	& $3701$ \\ \hline
$\beta_2$	& $(T^4+7T^3+15T^2+14T+9)/2$		& $-3701$ \\ \hline
$\beta_3$	& $(14T^5+155T^4+497T^3+439T^2-174T+143)/2$	& $3701^3$	\\ \hline
\end{tabular}
\end{center}
\caption{Some elements of $L$.}
\label{elements}
\end{table}

Let $L_p = \Q_p[T]/(f(T))$ be the completion of $L$ at a prime $p$ of $\Q$.
This will be a field if and only if there is only one prime of $L$ above $p$, which happens when $p=2$, for instance.
For $p=3701$, we have
	$$L_{3701} \isom \Q_{3701} \times E \times F$$
where $E$ is a totally ramified extension of $\Q_{3701}$ of degree~2, and $F$ is the unramified extension of $\Q_{3701}$ of degree~3.
The element $T$ maps in $\Q_{3701}$ to something that is $1371$ modulo $3701$, and in $E$ to something that is $1727$ modulo the maximal ideal.

Finally we will need to know how~2 splits in the subfield $K$ of $M$ corresponding to the subgroup $G$ of $S_6$ of elements that stabilize the partition $\{\{1,2,3\},\{4,5,6\}\}$ of $\{1,2,3,4,5,6\}$ into two indistinguishable subsets.
Since the orbit of $\{\{1,2,3\},\allowbreak\{4,5,6\}\}$ under the action of $S_6$ consists of $\binom{6}{3}/2=10$ partitions, $[K:\Q]=(S_6:G)=10$.
Let $\alpha_1,\ldots,\alpha_6$ be the roots of $f(x)$ in $M$, which we consider as a subfield of $\C$.
Then $\alpha_1 \alpha_2 \alpha_3 + \alpha_4 \alpha_5 \alpha_6 \in K$, and its conjugates are similar sums corresponding to the other partitions.
We can construct numerically the degree~10 polynomial $h(x)$ whose roots are these sums, and since these sums are the conjugates of an algebraic integer, the coefficients are integers, and we find the polynomial exactly:
{\footnotesize
	$$h(x) = x^{10} + 22 x^9 + 53 x^8 + 654 x^7 + 2186 x^6 + 8976 x^5 + 38705 x^4 + 89560 x^3 + 244664 x^2 + 565728 x + 477968.$$}%
This polynomial is irreducible over $\Q$, and it follows that $K=\Q(\alpha_1 \alpha_2 \alpha_3 + \alpha_4 \alpha_5 \alpha_6)$.
Finally, the prime~2 factors in $K$ as $\pp^4 \qq^2$ where $\pp$ is of degree $1$ and $\qq$ is of degree $3$, so in particular $h(x)$ has no zeros in $\Q_2$.

\section{The 2-descent on $\calC$}
\label{descent}

{}From now on, $J$ will denote the Jacobian of the curve $\calC$.
We will compute the Mordell-Weil rank of $J$ by performing the $2$-descent outlined in Section~\ref{generalities}.
Since $f$ has discriminant $2^{12} \cdot 3701$, we take $S = \{2,3701,\infty\}$, which contains all possible primes of bad reduction for $J$.
(In fact, our curve has good reduction at~$2$, because substituting $y=2z+x^3+x+1$ and dividing by~$4$ yields the model
	$$z^2 + x^3 z + x z + z = 2 x^5 + 5 x^4 + 5 x^3 + x^2 + x,$$
which has bad reduction only at $3701$.
But because we are doing a $2$-descent, we must include $2$ in $S$ anyway.)
Let $J(\Q)\tors$ denote the torsion subgroup of the finitely generated abelian group $J(\Q)$.

\begin{prop}
\label{trivialtorsion}
$J(\Q)\tors$ is trivial.
\end{prop}

\begin{proof}
For any prime $p$ of good reduction for $J$, the reduction mod $p$ map from $\jq$ to $ J({\F}_p)$ is injective on torsion.
(See~\cite{katz}, for example.)
By~\cite[p.\ 822]{GG},
	$$\# J(\F_p) = \frac{1}{2} \# \calC (\F_{p^{2}}) + \frac{1}{2}( \# \calC(\F_{p})^{2}) - p.$$
This can be obtained alternatively by evaluating the characteristic polynomial at~$1$.  (For a formula for the characteristic polynomial, see the proof of Proposition~\ref{noendomorphisms} in Section~\ref{nonmodularity}.)
Using this, we find $\# \jft = 9$ and $\# J ({\F}_5 ) = 41$.
But $\gcd(9,41)=1$, so $\# J(\Q)\tors=1$.
\end{proof}

An immediate corollary is that $[\infty^+-\infty^-]$ generates an infinite subgroup of $J(\Q)$ so the rank of $J(\Q)$ is at least~$1$.
Also, the fact $\# J(\F_5) = 41$ easily implies the following:

\begin{prop}
$J$ is not isogenous over $\Q$ to a product of two elliptic curves $E_1,E_2$ over $\Q$.
\end{prop}

\begin{proof}
If $J$ were isogenous over $\Q$ to $E_1 \times E_2$, then $E_1$ and $E_2$ would have good reduction at 5 as well.
Also $\# J(\F_5) = \# E_1(\F_5) \# E_2(\F_5)$, so $\# E_1(\F_5)$ and $\# E_2(\F_5)$ would be~1 and 41 in some order.
Both of these violate Hasse's bound
	$$|\# E(\F_p) - (p+1)| \le 2 \sqrt{p}.$$
\end{proof}

In Proposition~\ref{noendomorphisms} of Section~\ref{nonmodularity}, we will prove the much stronger result that $J$ is absolutely simple, and that $J$ has no nontrivial endomorphisms over $\C$.
This rules out the possibility of reducing the computation of the rank of $J(\Q)$ to the computation of ranks of elliptic curves, so we will need to use the general method outlined in Section~\ref{generalities}.
We proceed by first calculating the groups $G,G',H,H'$ of Section~\ref{generalities} for our curve.

\begin{lemma}
\label{unramifiedgroup}
The images of the~8 elements listed in Table~\ref{elements} in $L^\ast/L^{\ast 2}$ are a basis for the $\F_2$-vector space $G$.
\end{lemma}

\begin{proof}
If $l \in L^\ast$, then the field extension $L(\sqrt{l})/L$ is unramified at all finite primes of $L$ except possibly those occuring in $l$ and those above~$2$.
Hence the images of the~8 elements in Table~\ref{elements} are in $G$.
On the other hand, if $l \in L^\ast$ maps to something in $L^\ast/L^{\ast 2}$ outside the span of these~8 elements, then there must be a prime of $L$ not above $2$, $3701$ or $\infty$ that occurs to an odd power in the prime factorization of $l$, and then $L(\sqrt{l})/L$ is ramified at that prime.
\end{proof}

\begin{lemma}
\label{hgroup}
The images of $u_1$ and $u_3 \beta_1 \beta_2$ in $L^\ast/L^{\ast 2}\Q^\ast$ form a basis for the $\F_2$-vector space $H$.
\end{lemma}

\begin{proof}
If a prime $p$ other than $2$ or $3701$ occurs to an odd power in the factorization of $q \in \Q^\ast$, then $\Q(\sqrt{q})/\Q$ is ramified at $p$, and $L/\Q$ is unramified at $p$, so $L(\sqrt{q})/L$ is ramified at any prime above $p$ and $q \not\in G$.
On the other hand, by Table~\ref{badprimes}, the images of $-1$, $2$, $3701$ equal the images of $-1$, $u_2$, and $\beta_1 \beta_3$ in $G$.
Therefore $G' \subset L^\ast/L^{\ast 2} \Q^\ast$ is the quotient of $G$ by the latter three elements, and the images of $u_1,u_3,\alpha,\beta_1,\beta_2$ form a basis.
By Table~\ref{elements}, the kernel of the norm map from $G'$ to $\Q^\ast/\Q^{\ast 2}$ is the subspace generated by $u_1$ and $u_3 \beta_1 \beta_3$.
\end{proof}

\begin{lemma}
\label{tablelemma}
The last three columns of Table~\ref{badprimes} are accurate.
\end{lemma}

\begin{proof}
The nontrivial 2-torsion points of $J$ over $\Qbar_p$ are of the form $[(\alpha_i,0)+(\alpha_j,0) -\infty^+ - \infty^-]$, where $\alpha_i,\alpha_j$ are two of the six zeros of $f(x)$.
Over $\Q_2$, $f(x)$ is irreducible, so $\Gal(\Qbar_2/\Q_2)$ acts transitively on the six zeros, and hence no pair can be Galois-stable.
Thus $J(\Q_2)[2]$ is trivial.

Over $\Q_{3701}$, $f(x)$ factors into polynomials of degrees $1,2,3$.
Here the only pair of zeros that is stable under $\Gal(\Qbar_{3701}/\Q_{3701})$ is the pair of zeros of the quadratic factor.
Hence $J(\Q_{3701})[2]$ has one nontrivial point.

Over $\R$, $f(x)$ factors into polynomials of degrees $1,1,2,2$.
The pairs of zeros stable under complex conjugation are the pair of real zeros, and the pairs of zeros of each quadratic factor.
Hence $J(\R)[2]$ has three nontrivial points.

The multiplication-by-2 map on $J(\Q_p)$ is an $n$-to-1 map onto its image, where $n=\#J(\Q_p)[2]$, and locally it multiplies Haar measure by $|2|_p^2$ since $J(\Q_p)$ is a 2-dimensional Lie group over $\Q_p$.
Hence the measure of $2J(\Q_p)$ is $|2|_p^2/n$ times the measure of $J(\Q_p)$, so
	$$\# J(\Q_p)/2J(\Q_p) = |2|_p^{-2} \cdot \#J(\Q_p)[2],$$
which gives the values of the second to last column of Table~\ref{badprimes}.

{}From the factorization of~2 in $L$, we know that $f(x)$ has no roots in $\Q_2$.
{}From Section~\ref{numberfield}, the polynomial $h(x)$ has no roots in $\Q_2$, so there is no $\Gal(\Qbar_2/\Q_2)$-stable partition of the roots of $f(x)$ into two indistinguishable 3-element subsets.
Thus by Proposition~\ref{halves}, $2J(\Q_2)$ has index~2 in the kernel of $x-T$ on $J(\Q_2)$.
On the other hand, $f(x)$ has a zero in $\Q_{3701}$ and in $\R$, so Proposition~\ref{halves} implies that the kernel of $x-T$ on $J(\Q_p)$ equals $2J(\Q_p)$ for $p=3701$ or $p=\infty$.
\end{proof}

\begin{table}
{\tiny
\begin{center}
\begin{tabular}{|c||c|c|c|c|c|}
$p$	& ${(e_i,f_i)}$	& Factorization in $L$	& $\#J(\Q_p)[2]$	& $\#J(\Q_p)/2J(\Q_p)$ 	& $\#J(\Q_p)/\ker(x-T)$	\\ \hline \hline
2	& $(2,3)$	& $\alpha^2 u_2$	& 1	& 4	& 2	\\ \hline
3701	& $(1,1);(2,1);(1,3)$	& $\beta_1 \beta_2^2 \beta_3$	& 2	& 2	& 2	\\ \hline
$\infty$& $(1,1);(1,1);(2,1);(2,1)$	&& 4	& 1	& 1	\\ \hline
\end{tabular}
\end{center}
}
\caption{The primes in $S$.}
\label{badprimes}
\end{table}

Next we will need to find generators for $J(\Q_p)/\ker(x-T)$ for each prime $p$ in $S$.

\begin{lemma}
\label{generators}
The 1-dimensional $\F_2$-vector spaces 
$$J(\Q_2)/\ker(x-T) \hbox{\quad and\quad}
J(\Q_{3701})/\ker(x-T)$$ are generated by 
$$[(2,\sqrt{881})-\infty^-] \in J(\Q_2) \hbox{\quad and\quad}
[(-4,\sqrt{185})-\infty^-] \in J(\Q_{3701}),
$$ respectively.
\end{lemma}

\begin{proof}
For $p=2$, we have $881 \equiv 1 \pmod{8}$, so Hensel's Lemma implies that $(2,\sqrt{881})$ is in $\calC(\Q_2)$.  (Fix a square root.)
Thus it will suffice to show that $2-T \not\in L_2^{\ast2} \Q_2^\ast$.
Let $g(x)$ be the characteristic polynomial of $2-T$.
PARI tells us that there is only one prime above 2 in the number field generated by a root of $g(x^2)$, and it follows that $L_2(\sqrt{2-T})$ is a field of degree~12, so $2-T \not\in L_2^{\ast2}$.
Similarly, for each $r \in \{\pm 1,\pm 2,\pm 3,\pm 6\}$, a set of representatives for $\Q_2^\ast/\Q_2^{\ast2}$, we can check that $r(2-T) \not\in L_2^{\ast2}$, and it follows that $2-T \not\in L_2^\ast/L_2^{\ast2} \Q_2^\ast$.
(It should be remarked here, that it took PARI a few hours to do these calculations with degree~12 number fields.
We speculate that this is because the PARI command {\verb+initalg+}, which must precede the command {\verb+primedec+} that computes the decomposition of primes, computes many other pieces of information that are irrelevant for our purposes.  Of course, there are other methods that could be used to test if an element $x$ of $L_2^\ast$ is a square; for instance, if $x$ is a unit, this is determined by $x \bmod 8$.)

For $p=3701$, we first verify that the Legendre symbol $\left(\frac{185}{3701}\right)$ is~1, so Hensel's Lemma implies that $(-4,\sqrt{185}) \in \calC(\Q_{3701})$.
To complete the proof, we must check that $-4-T \not\in L_{3701}^{\ast2} \Q_{3701}^\ast$.
This time we can avoid the PARI computations with degree~12 number fields by exploiting the decomposition of $L_{3701}$ into fields.
As before, it suffices to prove that $r(-4-T) \not\in L_{3701}^{\ast2}$ for $r \in \{1,2,3701,2 \cdot 3701\}$, which is a set of representatives for $\Q_{3701}^\ast/\Q_{3701}^{\ast2}$, since $\left(\frac{2}{3701}\right)=-1$.
Now $-4-T$ maps in $\Q_{3701}$ to something that is $-4-1371 = -1375$ modulo $3701$, and $\left(\frac{-1375}{3701}\right)=1$, so by Hensel's Lemma, $-4-T$ maps to a square in $\Q_{3701}$, and $r(-4-T)$ can map to a square in $\Q_{3701}$ only if $r=1$.
On the other hand $-4-T$ maps in $E$ (the ramified field component of $L_{3701}$) to something that is $-4-1727=-1731$ modulo the maximal ideal, and $\left(\frac{-1731}{3701}\right)=-1$, so $-4-T$ does not map to a square in $E$, and hence $-4-T \not\in L_{3701}^{\ast2}$.
Thus $-4-T \not\in L_{3701}^{\ast2} \Q_{3701}^\ast$, and we are done.
\end{proof}

\begin{lemma}
\label{trivial}
$J(\Q)/\ker(x-T)$ is trivial.
\end{lemma}

\begin{proof}
By Proposition~\ref{landsinh}, diagram~(\ref{cd}) and Lemma~\ref{generators}, $J(\Q)/\ker(x-T)$ maps into the subgroup $H'$ of $H$ that maps in $L_2^\ast/L_2^{\ast 2} \Q_2^\ast$ into the group generated by $2-T$, in $L_{3701}^\ast/L_{3701}^{\ast 2} \Q_{3701}^\ast$ into the group generated by $-4-T$, and in $L_\infty^\ast/L_\infty^{\ast 2} \R^\ast$ to the identity.
So it will suffice to show that $H'$ is trivial.

First of all, the $\beta_2$-adic valuation $E \rightarrow \Z$ induces a map $v: L_{3701}^\ast/L_{3701}^{\ast 2} \Q_{3701}^\ast \rightarrow \Z/2\Z$, since the ramification index of $E$ over $\Q_{3701}$ is~2.
By Section~\ref{numberfield}, $-4-T$ maps in $E$ to something that is $-1731$ modulo the maximal ideal, so $v$ is trivial on the image of $-4-T$.
But $v$ maps the two generators $u_1$ and $u_3 \beta_1 \beta_2$ of $H$ to $0$ and $1$, respectively, so $H'$ is contained in the image of $\{1,u_1\}$.

The same method used in the proof of Lemma~\ref{generators} to show that $2-T$ was nontrivial in $L_2^\ast/L_2^{\ast 2} \Q_2^\ast$ shows that $u_1$ and $u_1(2-T)$ are nontrivial there, so $u_1$ does not map into the subgroup of $L_2^\ast/L_2^{\ast 2} \Q_2^\ast$ generated by $2-T$.
Thus $H'$ is trivial.
(The information from the prime $\infty$ was not used, but in fact it would not have helped either, since the kernel of the norm from $L_\infty^\ast/L_\infty^{\ast 2} \R^\ast$ to $\R^\ast/\R^{\ast 2}$ is trivial.)
\end{proof}

\begin{thm}
\label{rankone}
$J(\Q) \isom \Z$ as an additive group.
\end{thm}

\begin{proof}
By Proposition~\ref{halves}, $2J(\Q)$ has index~$2$ in $\ker(x-T)$, so
by Lemma~\ref{trivial}, we have
$\# J(\Q)/2J(\Q) = 2$.
By Proposition~\ref{trivialtorsion}, $J(\Q) \isom \Z^r$ for some $r \ge 1$.
Then $J(\Q)/2J(\Q) \isom (\Z/2\Z)^r$, so by the above, $r=1$.
\end{proof}

\section{Applying Chabauty's method}
\label{chabauty}

We recall the following consequence of Chabauty's result~\cite{CHA},
which gives a way of deducing information about 
the $\Q$-rational points on a curve from its Jacobian.

\begin{prop}
\label{prop1}
Let $C$ be a curve of genus~$g$
defined over $\Q$, whose Jacobian has Mordell-Weil
rank $\leq g-1$. Then $C$ has only finitely many
$\Q$-rational points.\end{prop}

This is a weaker result than
Faltings' Theorem; however, when applicable, Cha\-bau\-ty's
method can often be used to give good bounds
for the number of points on a curve.
Recent work in
Coleman~\cite{COL} (see also~\cite{MCA,MCB}) has improved 
Chabauty's technique; however, the bounds obtained seem
only rarely to resolve $\cq$ completely. 
For our curve $\calC$, the best bound that can
be obtained from the results in~\cite{COL} is
that $\# \cq  \leq 9$.
We shall adopt a more flexible approach that will
allow us to sharpen this bound to 6, as required.
It is hoped that a generalisation of the following ideas to any
curve of genus~2 over a number field will at some stage
be presented in~\cite{FLF}, but we make no direct
use of this, and present a largely self contained account
tailored to the needs of our specific example.
We shall, however, need to refer to the equations
in~\cite{FCPS,FCRE} relating to the Jacobian and formal group.
We shall first establish a few easily computed facts about
$\jq$.
Let $D = [ \infty^+ - \infty^- ] \in J(\Q)$.

\begin{lemma}
\label{le:D=kE}
We have $\jq = \langle E \rangle$, for some $E \in J(\Q)$ of infinite order, and $D=k \cdot E$ with $3 \notdiv k$.
\end{lemma}

\begin{proof}
By Theorem~\ref{rankone}, we can pick a generator $E$ for $\jq \isom \Z$.
To complete the proof, we must show that $D \not\in 3\jq$.
Since $ \jft $ is a cyclic group of size~9
generated by $\td$, the reduction of $D$ mod~3,
we find that $\td \not\in 3\jft$, from which it follows
that $D \not\in 3\jq$,
as required. 
\end{proof}
 
It would be nice to have the theory of heights sufficiently well
developed to determine whether $k=\pm 1$, which would
give $\jq = \langle D \rangle$. However, the method
in~\cite{FLH}
would require significant enhancements before it could
realistically be applied to $\calC$.
In fact, all of our local arguments will be 3-adic and
so the fact that $D\not\in 3\jq$ will turn out to
be sufficient for our purposes.

Table~\ref{divisormultiples} lists the first 11 multiples of $D$, which
will be relevant to our later computations. The last column
gives the corresponding multiples of
$\td$, the reduction of $D$ mod 3 to $\jft$.
For simplicity, we represent multiples of $D$ in ${\rm Pic}^{2}(\calC)$,
where ${\rm Pic}^{0}(\calC)\cong {\rm Pic}^{2}(\calC)$ by $[V]\mapsto [V+
\infty^{+}+\infty^{-}]$. In abuse of notation we will write $D=[\infty^{+}
-\infty^{-}]$ in ${\rm Pic}^{0}(\calC)$ and $D=[\infty^{+}+\infty^{+}]$
in ${\rm Pic}^{2}(\calC)$.
In the table, $P = (-2+\frac{1}{3}\sqrt{33}, 
               -\frac{17}{3}+\frac{10}{9}\sqrt{33} )$ and
$Q = ( -\frac{1}{2} + \frac{1}{6}\sqrt{-87},
        \frac{22}{3} + \frac{5}{9}\sqrt{-87} )$, and $\overline{P}$ and $\overline{Q}$ are their algebraic conjugates.

\begin{table}
\begin{center}
\begin{tabular}{|c||c|c|c|c|c|}
$n$	& $n \cdot D$	& $n \cdot \td$	\\ \hline \hline
0	& $\calO$	& $\calO$		\\ \hline
1	& $[\infty^+ + \infty^+]$	& $[ \infty^+ + \infty^+]$ \\ \hline
2	& $[ (0,1) + (-3,1) ]$		& $[ (0,1) + (0,1) ]$ \\ \hline
3	& $[ (0,-1) + \infty^- ]$	& $[ (0,-1) + \infty^- ]$ \\ \hline
4	& $[ (0,-1) + \infty^+ ]$	& $[ (0,-1) + \infty^+ ] $ \\ \hline
5	& $[ (-3,1) + \infty^- ]$	& $[ (0,1) + \infty^- ]$ \\ \hline
6	& $[ (-3,1) + \infty^+ ]$	& $[ (0,1) + \infty^+ ] $ \\ \hline
7	& $[ (0,-1) + (0,-1) ]$		& $[ (0,-1) + (0,-1) ] $ \\ \hline
8	& $[ P+ \overline{P}]$		& $[ \infty^- + \infty^- ]$ \\ \hline
9	& $[ (0,-1) + (-3,1) ]$		& $\calO$ \\ \hline
10	& $[ Q + \overline{Q}]$		& $[ \infty^+ + \infty^+] $ \\ \hline
11	& $[ (-3,1) + (-3,1) ]$		& $[ (0,1) + (0,1) ]$ \\ \hline
\end{tabular}
\end{center}
\caption{The first~11 multiples of $D$ and $\td$.}
\label{divisormultiples}
\end{table}

The multiples $\ell \cdot D$, for $\ell = -1,\ldots, -11$,
can be deduced from the above by using the rule
that $-[ (x_1,y_1) + (x_2,y_2)] = [ (x_1,-y_1)+ (x_2, -y_2)]$.
The divisor $9\cdot D$, which is in the kernel of reduction
mod 3, will play a special role, and so we denote:
	$$D' = 9\cdot D = [ (0,-1) + (-3,1) ].$$
\noindent The following lemma is immediate from the
fact that the $k$ of Lemma~\ref{le:D=kE} is coprime to 3.

\begin{lemma}
\label{le:lD+ME'}
Let $E$ be as in Lemma~\ref{le:D=kE}, and let $E' = 9\cdot E$.
Then any member of $\jq$ can be written uniquely
as $\ell \cdot D + m\cdot E'$, for some $\ell , m \in \Z$,
with $-4\leq \ell \leq 4$. \end{lemma}
 
If we now let:
	$$\calMt = \hbox{ the kernel of the reduction map from } 
			\jqt \hbox{ to } \jft ,$$
\noindent
then $\calMt$ contains no $k$-torsion, since $3\notdiv k$, and
there is a well defined map $1/k$ on $\calMt$ that takes any
$D_0\in \calMt$ to the unique $E_0\in \calMt$ such that $D_0 = k\cdot E_0$.
We can therefore legitimately say that any divisor in
$\jq$ can be written in the form:
\begin{equation}
\label{eqnf}
	\ell\cdot D + n\cdot D',\hbox{ with }-4\leq \ell \leq 4,
	\ n = m/k,\ 3\notdiv k, 
\end{equation}
\noindent where it is to be understood that $1/k$ refers
to the above 3-adic map on $\calMt$. Here, $n$ need not be
a rational integer, but must still be a $3$-adic
integer, which will be sufficient for our purposes. 

Our next observation is that 
$\cq$ is
in 1-1 correspondence with 
the members of $\jq$
that have the special form: $[ P+ P]$.
{}From Table~\ref{divisormultiples},
we see that all of the known $\Q$-rational points correspond
to: $\pm D$, $\pm 7\cdot D$ and $\pm 11 \cdot D$.
Suppose now that we have a divisor $D_0 \in \jq$ that
is of the special form $[ P+ P]$; we can write $D_0
= \ell \cdot D + n \cdot D'$ as in equation~(\ref{eqnf}).
If $D_0$ were in $\calMt$
(that is, $\ell = 0$),
then $\widetilde P$ would have to be of the form
$(x,0)$, which is impossible since the sextic $f(x)$ has
no roots in $\ft$. 
Otherwise, the reduction of $D_0$,
which is also the reduction of $\ell D$,  must be of the form
$[ \widetilde P + \widetilde P]$, giving that $ \pm 1 , \pm 2$
are the only possibilities for $\ell$.
Suppose we can show that $D + n\cdot D'$ is of the form $[ P+P]$
only when $n=0$ and that $2\cdot D + n\cdot D'$
is of that form only when $n = \pm 1$.
Using the fact that $-[ (x,y) + (x,y) ] = [ (x,-y)+ (x,-y) ]$,
it would then follow that $-D + n\cdot D'$ is of that form
only when $n=0$ and that $-2\cdot D + n\cdot D'$
is of that form only when $n=\pm 1$. This would show that
$\cq$ consists only of the 6 known points.
We summarise the above in the following lemma.

\begin{lemma}
\label{C(Q)=6}
Let $\calMt$ be the kernel of the reduction map from
$\jqt$ to $\jft$.
Let $D_1 = D = [ \infty^+ + \infty^+ ]$
and $D_2 = 2\cdot D = [ (0,1)+ (-3,1)]$. Then $D' = 9\cdot D =
[ (0,-1)+ (-3,1) ] \in \calMt$. Suppose that, for
$n = m/k,\ m, k \in \Z , 3\notdiv k$,
we have $D_1 + n\cdot D'$ of the form $[ P+P]$ only when
$n=0$, and $D_2 + n\cdot D'$ is of that form only when
$n = \pm 1$. Then $\# \cq = 6$. \end{lemma}

For each $D_i$, $i=1,2$, our strategy will be to derive, to a sufficient
degree of 3-adic accuracy,
a power series $\theta_i(n) \in \zt [[n]]$ that must be
satisfied by $n$ whenever $D_i + n\cdot D'$ is of the
form $[ P+P]$. We shall show the stronger result that
the known solutions to $\theta_i(N)$ give all of the
solutions $n\in \zt$. 
The following standard theorem of Strassman
is proved in~\cite[p.62]{CA1}.

\begin{thm}
\label{Strassman}
Let $\theta (X) = c_0 + c_1 X + c_2 X^2 + \ldots
\in \zp [[X]] $ satisfy $c_n \rightarrow 0$ in $\qp$.
Define $r$ uniquely by: $\vert c_r \vert _p \geq\vert c_i
\vert_p$ for all $i\geq0$, and $\vert c_r \vert_p
> \vert c_i
\vert_p$ for all $i>r$.
Then there are at most $r$ values of $x\in \zp$ such that
$\theta (x) = 0$.\end{thm}

In order to derive the power series $\theta_1(n)$ and $\theta_2(n)$,
we shall make use of the formal group.
As remarked in Section~\ref{hyperelliptic}, $\calC$ cannot be put in the simpler form $y^2 = (\text{\sl quintic in } x)$, so instead of using the
development of the formal group in~\cite{GRA},
we must use the general $y^2 = (\text{\sl sextic in } x)$\ 
development as in~\cite{FCPS,FCRE}. The derivation of the equations that we shall use for both the formal group law and the global group law
are described in~\cite{FCRE}. These equations for a general curve of genus~2, are available
at: \verb+ftp.liv.ac.uk+ in the directory
$\;\lower3.7pt\hbox{$\widetilde{ }$}\hskip3pt $\verb+ftp/pub/genus2+
by anonymous ftp. 
First note that for any curve of genus~2
\begin{equation}
\label{eqng}
y^2 = f_6 x^6 + f_5 x^5 + f_4 x^4 + f_3 x^3 + f_2 x^2
+ f_1 x + f_0, \ \ f_i \in \Z ,
\end{equation}
the following functions $s_1,s_2$ of a point $D_0 = [ (x_1,y_1)+ (x_2,y_2)]\in \jq$ can be used as a pair of local parameters at $\calO$:
\begin{align}
\label{eqnh}
s_1 &= (G_1 (x_1 ,x_2 )y_1 - G_1 (x_2 ,x_1 )y_2 )(x_1 - x_2 ) 
/(F_0 (x_1 ,x_2 )-2y_1 y_2 )^2, \\
s_2 &= (G_0 (x_1 ,x_2 )y_1 - G_0 (x_2 ,x_1 )y_2 )(x_1 - x_2 )
/(F_0 (x_1 ,x_2 )-2y_1 y_2 )^2,
\end{align}
\noindent where
\[
\begin{array}{rl}
F_0 (x_1 ,x_2 )= & 2 f_0 + f_1 (x_1 +x_2 ) + 2 f_2 (x_1 x_2 )
+ f_3 (x_1 x_2) (x_1 +x_2 ) \\
 & + 2 f_4 (x_1 x_2 )^2 + f_5 (x_1 x_2 )^2 (x_1 + x_2 )
+2 f_6 (x_1 x_2 )^3, \\
G_0 (x_1 ,x_2 )= & 4 f_0 + f_1 (x_1+3 x_2)+f_2(2x_1 x_2+2 x_2 ^2 )
   + f_3(3x_1 x_2 ^2 + x_2 ^3 )  \\
 & +4 f_4(x_1 x_2 ^3)+f_5(x_1 ^2 x_2 ^3 + 3 x_1 x_2 ^4)
  + f_6(2 x_1 ^2 x_2 ^4 + 2 x_1 x_2 ^5 ), \\
G_1 (x_1 ,x_2 )= & f_0(2x_1+2x_2)+f_1(3x_1 x_2 + x_2 ^2) + 4f_2(x_1 x_2 ^2)
  + f_3(x_1 ^2 x_2 ^2 +3x_1 x_2^3) \\
 & +f_4(2x_1^2 x_2^3 +2x_1 x_2^4) +f_5(3x_1^2 x_2^4 +x_1 x_2^5)
+4f_6(x_1^2 x_2^5).
\end{array}\]
\medskip

The following lemma summarises the information we need from
\cite{FCPS,FCRE} and introduces the standard formal exponential
and logarithm maps on the formal group.

\begin{thm}
\label{formal}
Let $C$ be as in~(\ref{eqng}).
There is a formal group law with respect to the local
parameters of equation~(\ref{eqnh}), given by
$\calF =\begin{pmatrix}
		\calF_1 \\ \calF_2
	\end{pmatrix}$
where $\calF_1 , \calF_2$ are power series
in $s_1,s_2,t_1,t_2$ defined over
$\Z$, which contain terms only of odd degree.
Define the {\it formal exponential of}
$\calF$ as $E = \begin{pmatrix} E_1 \\ E_2 \end{pmatrix}$,
where $E_1$, $E_2$ are power series in {\bf s}
over $\Q$, by:
$E({\bf s}) = {\bf s} + \hbox{ terms of higher degree}$,
and $E({\bf s} + {\bf t}) = {\calF}(E({\bf s}), E({\bf t}))$.
Similarly define the {\it formal logarithm of} $\calF$ as 
$L = \begin{pmatrix} L_1 \\ L_2 \end{pmatrix}$
where $L_1$, $L_2$ are power series in {\bf s}
over $\Q$,  by: $L(E({\bf s})) = {\bf s}$,
or
equivalently:
$L({\bf s}) = {\bf s} + \hbox{ terms of higher degree}$,
and $L({\calF}({\bf s},{\bf t})) = L({\bf s}) + L({\bf t})$.
Then each of $E_1$,$E_2$,$L_1$,$L_2$
can be written in the form: $\sum (a_{ij}/i!j!) s_1^i s_2^j$,
where $a_{ij} \in \Z$ and $a_{ij} = 0$ when $i+j$ is even.
Let $p$ be a prime of good reduction, and let
$A,B,C$ be in $\calMp$, the kernel of reduction from
$\jqp$ to $\jfp$, with
$C = A + B$.
Suppose now that ${\bf s} = \begin{pmatrix} s_1 \\ s_2 \end{pmatrix}$ are
the local parameters corresponding to $A$, and
similarly {\bf t}, {\bf u} those for $B$, $C$ respectively.
Then each $s_i,t_i,u_i \in p\zp$ and 
${\calF}({\bf s},{\bf t})$ converges in $p\zp$ with
${\bf u} = {\calF}({\bf s},{\bf t})$.\end{thm}

The power series $\calF$ gives a description of the  
group law on $\calMp$. 
It is described in~\cite{FCPS,FCRE} how to compute
terms of the formal group up to terms of arbitrary degree.
We require here the formal group up to terms of degree~3 in~{\bf s}:
\begin{align*}
{\calF}_1 &= s_1 + t_1 + 2f_4s_1^2t_1 + 2f_4s_1t_1^2
- f_1s_2^2t_2 - f_1s_2t_2^2 + (\hbox{degree }\geq5) \\
{\calF}_2 &= s_2 + t_2 + 2f_2s_2^2t_2 + 2f_2s_2t_2^2
- f_5s_1^2t_1 - f_5s_1t_1^2 + (\hbox{degree }\geq5)
\end{align*}
For any $A$ in
$\calMp$, with local
parameter ${\bf s}$, note also that the power series $E({\bf s})$
and $L({\bf s})$ converge in $p\zp$ also, since 
$\vert s_1 \vert _p , \vert s_2 \vert _p \leq p^{-1}$ and so
$\vert s_1^i s_2^j / i! j! \vert_p$ converges to $0$
as $i+j \rightarrow \infty$. 
Once terms of the formal group have been computed, the terms
of $E$ and $L$ may be computed inductively from their definitions.
We shall again only require terms up to degree~3 in~{\bf s}:
$$ \begin{array}{lll}
L_1({\bf s}) = s_1 + \frac{1}{3}(-2f_4 s_1^3 + f_1 s_2^3) +\ldots
& \;\;\;\;\;\;\;\; &
E_1({\bf s}) = s_1 + \frac{1}{3}( 2f_4 s_1^3 - f_1 s_2^3) +\ldots  \\
L_2({\bf s}) = s_2 + \frac{1}{3}(-2f_2 s_2^3 + f_5 s_1^3) +\ldots 
& \;\;\;\;\;\;\;\; &
E_2({\bf s}) = s_2 + \frac{1}{3}( 2f_2 s_2^3 - f_5 s_1^3) +\ldots 
\end{array}
$$
Let us now return to our specific curve $\calC$ of equation~(\ref{curveequation}).
The local parameters of $D' = [ (0,-1) + (-3,1) ] \in \calMt$
are determined by substituting $x_1=0,y_1=-1,x_2=-3,y_2=1$
into~(\ref{eqnh}), giving: $s_1 = -9/14$ and $s_2 = 426/49$,
both of which have $3$-adic valuation less than or equal to $3^{-1}$. 
It is immediate that $L_1$,$L_2$ evaluated at $s_1 = -9/14$, $s_2 = 426/49$,
are both $3$-adic integers, and that (even after taking
denominators into account) the terms up to degree~3 determine
$L_1$,$L_2$ mod $3^4$. This gives: $L_1 \equiv 36 \ ( \hbox{mod } 3^4)$
and $L_2 \equiv 3 \ ( \hbox{mod } 3^4)$.
{}From the properties of $E$ and $L$ we see that $E(n\cdot L({\bf s}) )$
gives the local parameters $t_1,t_2$ for $T = n\cdot D'\in \calMt$,
where $n$ is as in~(\ref{eqnf}), and so is in $\zt$. This expresses
each of $t_1,t_2$ as members of $\zt [[n]]$, 
given (mod $3^4$) by:
\begin{equation}
\label{eqnk}
 t_1 \equiv 36 n + 27 n^3 \hbox{ and }
   t_2 \equiv 3n + 9n^3 \ \ (\hbox{mod } 3^4).
\end{equation}

Since any member of $\calMt$ is uniquely determined by its
local parameters, this describes $T = n\cdot D'$ as
a power series in $n$.
We now wish to describe $D_1 + T$ and $D_2 + T$,
where $D_1$,$D_2$ are as specified in Lemma~\ref{C(Q)=6}. 
Applying the standard global group law
to the sum $[ (x_1,y_1)+(x_2,y_2)] = D_1 + T$ 
gives (as described in~\cite{FCRE}) expressions
for $k_1,k_2,k_3\in \Z [[t_1,t_2]]$ such that
the triple $(k_1,k_2,k_3)$ is the same projectively
as $(1,x_1+x_2,x_1 x_2)$. 
The terms up to degree~3 in {\bf t} are:
$$
\begin{array}{l}
k_1 = 
-12 t_2-12 t_1^2+8 t_1 t_2+36 t_2^2
+8 t_1^3-72 t_1^2 t_2-48 t_1 t_2^2-8 t_2^3 +\cdots \\
k_2 = 
12 t_1+48 t_2-8 t_1^2-104 t_1 t_2-132 t_2^2
+72 t_1^3+648 t_1^2 t_2+408 t_1 t_2^2+104 t_2^3 +\cdots \\
k_3 = 
-6+4 t_1-72 t_1^2-24 t_1 t_2-4 t_2^2
-24 t_1^3-104 t_1^2 t_2-104 t_1 t_2^2-24 t_2^3 + \cdots\\
\end{array}
$$
\noindent On substituting~(\ref{eqnk}) into these expressions
gives each of $k_1,k_2,k_3$ as members of $\zt [[n]]$.
Now note that if a divisor $[ (x_1,y_1)+ (x_2,y_2) ]$  
is of the form $[ P+P]$ then 
$\theta_1(n) = k_2^2 - 4k_1k_3 = 0$. This gives:
	$$\theta_1(n) \in \zt [[n]], \hbox{ with } \theta_1(n) \equiv 27n \ (\hbox{mod } 3^4),$$
\noindent
where $\theta_1(n) = 0$ if $D_1 + n\cdot D'$ is of the
form $[ P+P]$.

Repeating the same process for $D_2$ first gives: 
{\small$$
\begin{array}{l}
k_1 = 
-2-12 t_1-40 t_2-16 t_1^2+64 t_1 t_2+100 t_2^2
-64 t_1^3-472 t_1^2 t_2-64 t_1 t_2^2-64 t_2^3 + \ldots \\
k_2 = 
6+36 t_1+116 t_2+52 t_1^2-224 t_1 t_2-392 t_2^2
+208 t_1^3+1408 t_1^2 t_2+72 t_1 t_2^2+160 t_2^3 + \ldots \\
k_3 = 
4 t_1+12 t_2+28 t_1^2+176 t_1 t_2+272 t_2^2
+32 t_1^3+208 t_1^2 t_2+104 t_1 t_2^2+16 t_2^3 + \ldots \end{array}
$$}%
which then gives:
$$\theta_2(n) \in \zt [[n]], \hbox{ with }  \theta_2(n) \equiv 36 + 27n + 18n^2 + 54n^3 + 27n^4 \ (\hbox{mod } 3^4),$$
\noindent
where $\theta_2 (n) = 0$ if $D_2 + n\cdot D'$ is of the  
form $[ P+P]$.
We are now in a position to prove the desired result.

\begin{thm}
\label{C(Q)}
The curve $\calC$ of equation~(\ref{curveequation}) has only the six $\Q$-rational points $(0,1)$, $(0,-1)$, $(-3,1)$, $(-3,-1)$, $\infty^+$, and $\infty^-$ listed in Table~\ref{cvalues}.
\end{thm}

\begin{proof}
The coefficient of $n$ in $\theta_1(n)\in \zt [[n]]$
has $3$-adic valuation strictly larger than all of the
other coefficients, and so by Strassman's Theorem
(Theorem~\ref{Strassman}) there is at most~1 solution, which
is the known solution: $n=0$.
For $\theta_2(n)$, we further reduce mod~$3^3$,
giving: $\theta_2(n) \equiv 9 + 18n^2$. By Strassman's
Theorem, there are at most 2 solutions, which must
be the 2 known solutions: $n=-1,-2$.
The result now follows from Lemma~\ref{C(Q)=6}.
\end{proof}

\section{Non-modularity of $C_0(5)$ and $C_1(5)$}
\label{nonmodularity}

Recall that $C_1(4)$ turned out to be isomorphic over $\Q$ to the modular curve $X_1(16)$.
Morton~\cite{morton4} asked whether $C_1(N)$ could be parameterized by modular functions also for $N>4$.
If $C_0(5)$ or $C_1(5)$ were isomorphic over $\C$ to $X_1(N)$ or $X_0(N)$, then $N$ could not be a multiple of $3701$, because by~\cite[Corollary 9.11]{knapp} the genus of $X_0(3701)$ already is $(3701-5)/12 = 308$, whereas by Table~\ref{genus}, $C_0(5)$ and $C_1(5)$ have genus~$2$ and $14$, respectively.
Hence $C_0(5)$ or $C_1(5)$ would have potential good reduction at $3701$.
Using Lange's theorem~\cite{lange} that potential good reduction of a geometrically connected smooth projective curve is inherited by any other such curve it surjects onto (or the more general result mentioned in the Appendix by Matignon and Youssefi to~\cite{youssefi} that the same is true for good reduction), we find that in either case, $C_0(5)$ would have potential good reduction at $3701$.
But it can be shown that this contradicts the fact that the exponent of $3701$ in the discriminant of $f(x)$ is~$1$, so neither $C_0(5)$ nor $C_1(5)$ is isomorphic over $\C$ to $X_0(N)$ or $X_1(N)$ for any $N \ge 1$.

We have been slightly sketchy in the previous argument, because below we will provide a complete proof for the stronger result that there is no surjective morphism from $X_1(N)$ to $C_0(5)$ or $C_1(5)$ for any $N \ge 1$, even over $\C$.
As before, let $J$ denote the Jacobian of $\calC=C_0(5)$.
Let $\End J$ denote the ring of endomorphisms of $J$ defined over $\C$.

\begin{prop}
\label{noendomorphisms}
$J$ is absolutely simple, and $\End J \isom \Z$.
\end{prop}

\begin{proof}
We will model our argument on that used in~\cite[Appendix A]{pyle}.
Suppose $p$ is a prime of good reduction for $J$.
Then reduction modulo $p$ embeds $\End J$ in $\End_{\Fbar_p} J$, the endomorphisms defined over $\Fbar_p$ of the reduced abelian variety over $\F_p$ (which we will also denote $J$).
By~\cite[Lemma~3]{merrimansmart}, the characteristic polynomial of the Frobenius endomorphism $\pi_p$ on $J$ is
\begin{equation}
\label{charpoly}
	X^4 - t X^3 + s X^2 - p t X + p^2,
\end{equation}
where
	$$t = p+1-\# \calC(\F_p), \qquad s = \frac{1}{2} \left[\# \calC(\F_p)^2+\# \calC(\F_{p^2}) \right] +p-(p+1)\# \calC(\F_p).$$
Moreover, it follows from~\cite[Theorem 8]{wm} that if the characteristic polynomial of $\pi_p^n$ is irreducible over $\Q$ for all $n \ge 1$, then $(\End_{\Fbar_p} J) \otimes \Q = \Q(\pi_p)$ is a number field of degree~$4$.

For $p=3$,~(\ref{charpoly}) becomes $X^4-X^2+9$, so the characteristic polynomial of $\pi_3^2$ is $(X^2-X+9)^2$.
Hence we move on to $p=5$, for which~(\ref{charpoly}) is $P(x)=X^4+X^3+9X^2+5X+25$.
This is irreducible over $\Q$, so $\Q(\pi_5) \isom \Q[X]/(P(x))$ is a number field of degree~$4$.
We wish to show that no positive power of $\pi_5$ lies in a proper subfield.
PARI tells us that the Galois group of $P(X)$ is dihedral of order $8$, so $\Q(\pi_5)$ has an automorphism $\sigma$ of order~$2$, even though it is not Galois over $\Q$.
By Galois theory, the (quadratic) fixed field $F$ of $\sigma$ is the only nontrivial subfield of $\Q(\pi_5)$.
We find that $\pi_5+\sigma(\pi_5)$ is a root of $x^2+x=1$, so $F=\Q(\sqrt{5})$.
If $\pi_5^n \in F$, then $\sigma(\pi_5)/\pi_5$ would be an $n$-th root of unity.
But PARI shows that the only roots of unity in $\Q(\pi_5)$ are $1$ and $-1$, and that $\sigma(\pi_5)/\pi_5$ is neither of these.
Thus we now know that $(\End_{\Fbar_5} J) \otimes \Q \isom \Q[X]/(P(X))$, which already is enough to imply that $J$ is absolutely simple.

The characteristic polynomial of $\pi_7$ is $R(X)=X^4+2X^3+4X^2+14X+49$, and exactly the same argument as in the previous paragraph shows $(\End_{\Fbar_7} J) \otimes \Q \isom \Q[X]/(R(X))$.
Now $(\End J) \otimes \Q$ embeds into both number fields $\Q[X]/(P(X))$ and $\Q[X]/(R(X))$, but PARI tells us that the only nontrivial subfield $F=\Q(\sqrt{5})$ of $\Q[X]/(P(X))$ is not a subfield of $\Q[X]/(R(X))$, so $(\End J) \otimes \Q = \Q$.
Thus $\End J=\Z$.
\end{proof}

Let $J_1(N)$ denote the Jacobian of $X_1(N)$.
We will write $\End_\Q A$ for the ring of endomorphisms defined over $\Q$ of an abelian variety $A$ over $\Q$.

\begin{prop}
\label{quotientsofj1}
Let $B$ be an absolutely simple abelian variety over $\C$ which is a quotient of $J_1(N)$ over $\C$.
Then the rank of $\End B$ over $\Z$ is $\dim B$ or $2 \dim B$.
\end{prop}

\begin{proof}
Let $A$ be a simple abelian variety over $\Q$ which is a quotient of $J_1(N)$ over $\Q$, and which contains $B$ in its decomposition into absolutely simple abelian varieties over $\C$ up to isogeny.
If $B$ is an elliptic curve with complex multiplication, the result is trivial, so assume this does not hold.
Then by~\cite[Theorem~$1$]{ribet}, $(\End A)\otimes \Q$ is a matrix algebra $D=\M_n(H)$ over a division algebra $H$ finite dimensional over its center $F$, and $E \defequal (\End_\Q A)\otimes \Q$ is a maximal subfield of $D$.
Moreover $[E:\Q]=\dim A$, and $[H:F]=r^2$ with $r=1 \text{ or } 2$.
Let $f=[F:\Q]$.
Since $E$ is a maximal subfield of $D$, $[E:F]=\sqrt{[D:F]}=\sqrt{n^2 r^2}=nr$, so
	$$\dim B = (\dim A)/n = [E:\Q]/n = [E:F] f/n = rf.$$
Finally,
	$$\rank(\End B)=[(\End B)\otimes\Q:\Q]=[H:\Q]=r^2 f$$
so $\rank(\End B)=r \dim B$, and we are done.
\end{proof}

\begin{thm}
\label{nomorphism}
Let $N \ge 1$.
There is no nonzero morphism of abelian varieties over $\C$ from $J_1(N)$ to $J$.
Thus there is no surjective morphism of curves from $X_1(N)$ to $C_0(5)$ or $C_1(5)$.
\end{thm}

\begin{proof}
By Proposition~\ref{quotientsofj1}, any $2$-dimensional quotient of $J_1(N)$ must have an endomorphism ring larger than $\Z$.
Thus the first statement follows from Proposition~\ref{noendomorphisms}.
Since $C_1(5)$ maps to $C_0(5)$, and since surjective maps on curves induce surjective maps on their Jacobians, the final statement follows from the first.
\end{proof}

For the modular curves $X_0(N)$ and $X_1(N)$, the Manin-Drinfeld theorem states the divisor class of the difference of two cusps is a torsion element in the Jacobian.
It is natural to ask whether the same is true for $C_0(N)$ and $C_1(N)$, with cusps replaced by points with $c=\infty$.
(All of these points are rational, as follows from the ``$q$-expansions'' in~\cite{morton4}.)
For $N=4$, the result holds, simply because $C_1(4)$ is isomorphic to $X_1(16)$ and the points with $c=\infty$ correspond to cusps.
But the result fails for $N=5$, even for the quotient $C_0(5)$, since the divisor class of the difference of two of its rational points at $c=\infty$ is a nonzero element of $J(\Q)$, and hence is not torsion, by Proposition~\ref{trivialtorsion}.

\section{Rational points and cycles of period $6$}
\label{period6}

We conclude the paper with a few remarks about the next unsolved case, $N=6$.
The curve $C_0(6)$ is of genus~$4$ (see Table~\ref{genus}) and is birational to the curve given by the equation $\tau_6(x,c)=0$, where
{\footnotesize\begin{equation*}
\thickmuskip=
.5\thickmuskip
\medmuskip=
.5\medmuskip
\begin{split}
	\tau_6(x,c)	&=	(-384 c - 592 {c^2} - 256 {c^3}) + 
   \left( 448 + 416 c - 304 {c^2} - 256 {c^3} \right)  x + 
   \left( 196 + 552 c + 480 {c^2} + 256 {c^3} \right)  {x^2} \\
	& \quad	+
   \left( 140 - 136 c + 160 {c^2} + 256 {c^3} \right)  {x^3} + 
   \left( 175 + 16 c + 112 {c^2} \right)  {x^4} + 
   \left( 49 + 16 c + 144 {c^2} \right)  {x^5} \\
	& \quad + 
   \left( 14 + 8 c \right)  {x^6} + \left( 2 + 24 c \right)  {x^7} - 
   {x^8} + {x^9}.
\end{split}
\end{equation*}}%
(This is taken from~\cite{morton4}.)
Recall that $x$ is the trace of a 6-cycle for $g(z)=z^2+c$.

\begin{table}
\begin{center}
\begin{tabular}{|c||c|c|c|}
$(x,c)$		& Generator of $6$-cycle	& $K$	& Conductor	\\ \hline \hline
$(0,0)$		& $\zeta_9$			& $\Q(\zeta_9)$	& $9$	\\ \hline
$(-1,-2)$	& $\zeta_{13}+\zeta_{13}^{-1}$	& $\Q(\zeta_{13}+\zeta_{13}^{-1})$	& $13$	\\ \hline
$(1,-2)$	& $\zeta_{21}+\zeta_{21}^{-1}$	& $\Q(\zeta_{21}+\zeta_{21}^{-1})$	& $21$	\\ \hline
$(-3,-4)$	& $(\zeta_7^2+\zeta_7^{-2}) + \frac{1+\sqrt{5}}{2} (\zeta_7+\zeta_7^{-1})$	& $\Q(\zeta_7+\zeta_7^{-1},\sqrt{5})$	& $35$	\\ \hline
$(-7/2,-71/48)$	& $-1+\frac{\sqrt{33}}{12}$	& $\Q(\sqrt{33})$	& $33$	\\ \hline
\end{tabular}
\end{center}
\caption{The known affine rational points on $C_0(6)$.}
\label{sixcycles}
\end{table}

For each rational number $x=r/s$ with $|r|,|s| \le 100$, we checked the polynomial $\tau_6(x,c)$ in $c$ for rational roots.
We then did the same with $x$ and $c$ reversed.
This let us find all affine rational points on $\tau_6(x,c)=0$ having at least one coordinate with numerator and denominator bounded by~$100$ (in absolute value).
These are listed in Table~\ref{sixcycles}.
Because each of these points in fact has a coordinate with numerator and denominator bounded by $7$, it seems reasonable to expect that we have found all the affine rational points.
(There are also $5$ points at infinity on the nonsingular model, and these are all rational.)

Each affine point on $C_0(6)$ corresponds to a $\GalQ$-stable $6$-cycle, whose elements generate abelian extensions of $\Q$ of degree dividing $6$.
Table~\ref{sixcycles} lists an element of this cycle for each known point (in terms of a primitive $n$-th root of unity $\zeta_n$), and also gives the abelian extension $K$ of $\Q$ it generates, together with its conductor.
(It is straightforward to verify these using PARI.)
In particular, note that none of the cycles are defined pointwise over $\Q$.
Therefore, if we have truly found all affine rational points on $C_0(6)$, then there is no quadratic polynomial $g(z) \in \Q[z]$ with a periodic point of exact period~$6$.

\section*{Acknowledgements}

We thank Greg Call for helping us trace the history of the problem mentioned in the first paragraph, Noam Elkies for a comment that let us check that $\calC$ actually had good reduction at $2$, Qing Liu for referring us to the theorems on good reduction of curves mentioned in Section~\ref{nonmodularity}, Patrick Morton for sharing his preprints with us, Michel Olivier for verifiying our number field computations unconditionally using a yet to be released version of PARI, Ken Ribet for suggesting to us that the implication $(\End J=\Z) \implies (J \text{ is not a modular quotient})$ in Section~\ref{nonmodularity} should follow easily from the results in~\cite{ribet}, and Michael Zieve for introducing us to the problems considered in this paper.


\end{document}